\documentclass[journal]{ieeetran}
\usepackage{cite}
\usepackage{graphicx}
\usepackage{epstopdf}
\usepackage{amsmath}
\usepackage{cases}
\usepackage[caption=false,font=footnotesize]{subfig}
\usepackage{fixltx2e}
\usepackage{amssymb}
\usepackage{mathtools}
\usepackage{bm}
\usepackage{multirow}
\usepackage{color}
\usepackage{algorithm}
\usepackage{algorithmic}
\usepackage{amsthm}
\usepackage{mathrsfs}

\newtheorem{assumption}{Assumption}
\newtheorem{proposition}{Proposition}
\newtheorem{remark}{Remark}

\newcommand{\tabincell}[2]{\renewcommand\arraystretch{0.9}\begin{tabular}{@{}#1@{}}#2\end{tabular}}

\begin{document}

\title{Nonintrusive Uncertainty Quantification of Dynamic Power Systems Subject to Stochastic Excitations}

\author{Yiwei~Qiu,~\IEEEmembership{Member,~IEEE},
        Jin~Lin,~\IEEEmembership{Member,~IEEE},
        Xiaoshuang~Chen, \\
        Feng~Liu,~\IEEEmembership{Senior Member,~IEEE},  and
        Yonghua~Song,~\IEEEmembership{Fellow,~IEEE}%

\thanks{Financial support came from the National Key R\&D Program of China (2018YFB0905200), the National Natural Science Foundation of China (51907099, 51761135015, 51577096), and the China Postdoctoral Science Foundation (2019M650676). (Corresponding author: Jin Lin)}
\thanks{Y. Qiu, J. Lin, X. Chen, and F. Liu are with the State Key Laboratory of Control and Simulation of Power Systems and Generation Equipment, Department of Electrical Engineering, Tsinghua University, Beijing, 100087, China. (email: 
 linjin@tsinghua.edu.cn)
}%
\thanks{Y. Song is with the Department of Electrical and Computer Engineering, University of Macau, Macau 999078, China, and the Department of Electrical Engineering, Tsinghua University, Beijing 100087, China. 
}

}

\maketitle

\begin{abstract}
  Continuous-time random disturbances (also called stochastic excitations) due to increasing renewable generation have an increasing impact on power system dynamics; However, except from the slow Monte Carlo simulation, most existing methods for quantifying this impact are \emph{intrusive}, meaning they are not based on commercial simulation software and hence are difficult to use for power utility companies. To fill this gap, this paper proposes an efficient and \emph{nonintrusive} method for quantifying uncertainty in dynamic power systems subject to stochastic excitations. First, the Gaussian or non-Gaussian stochastic excitations are modeled with an It\^{o} process as stochastic differential equations. Then, the It\^{o} process is spectrally represented by independent Gaussian random parameters, which enables the polynomial chaos expansion (PCE) of the system dynamic response to be calculated via an adaptive sparse probabilistic collocation method. Finally, the probability distribution and the high-order moments of the system dynamic response and performance index are accurately and efficiently quantified. The proposed nonintrusive method is based on commercial simulation software such as PSS/E with carefully designed input signals, which ensures ease of use for power utility companies. The proposed method is validated via case studies of IEEE 39-bus and 118-bus test systems.
\end{abstract}


\section{Introduction}
\label{sec:intro}

\IEEEPARstart{I}{n} recent years, the increasing penetration of renewable generation has posed increasing continuous-time disturbances on power systems. We call these \emph{stochastic excitations} in this paper. The stochastic excitations have a significant impact on system dynamic performance that needs to be quantified to provide information for system operation \cite{wang2015long,milano2013systematic}. Up to now, many works treat renewable generation as static random parameters that do not vary with time, but such simplification may lead to inaccuracy. This is because in these methods, uncertainty is only considered at the beginning of the dynamic process. Power system dynamics itself is, however, still modeled as deterministic, and stochastic variations during the process are neglected. On the other hand, researches on power system dynamic uncertainty quantification that consider continuous-time stochastic excitations are still in progress.

An adequate method for quantifying the uncertainty of a dynamic power system subject to stochastic excitations should satisfy the following three requirements: 1) \emph{Accuracy:} the result should reflect the probability distribution and high-order moments of the non-Gaussian uncertainty as well as the nonlinearity of the system. 2) \emph{Computational efficiency:} a fast method enables online evaluation of fast-changing operating conditions. 3) \emph{Ease of use for utility companies:} a nonintrusive method implemented in commercial simulation software is easier to use and more reliable than an intrusive one that uses alternative code to perform a dynamic analysis. Here, \emph{nonintrusive} means that an existing power system simulation tool can be directly employed without rewriting built-in models or numerical solvers; \emph{intrusive} is in contrast to this.

Generally, the Monte Carlo simulation (MCs) is the most straightforward way to power system dynamic uncertainty quantification  \cite{dong2012numerical,milano2013systematic,ortega2018stochastic}. MCs allows nonintrusive assessment of uncertainty. However, when higher precision is required, its inefficiency caused by the large number of samplings hinders its online application.

Emerging methods based on stochastic differential equations (SDEs) \cite{Apostolopoulou2016assessment,ju2018analytical,shi2018analytical,li2019analytic,Wang2013the,Chen2018Stochastic} also have drawn significant attention. These methods use SDEs to model dynamic power systems subject to stochastic excitations and then use stochastic calculus to analyze them, and can be traced back to the work of Wang and Crow \cite{Wang2013the} using the Fokker-Planck equation to describe the evolving probability densities of system states.

However, many of the SDE-based methods rely on the assumption of Gaussian uncertainty \cite{Apostolopoulou2016assessment,ju2018analytical,shi2018analytical,li2019analytic}, which could be unrealistic in practice; some use partial differential equations, e.g., the Fokker-Planck equation, to depict the evolution of the probability distribution \cite{Wang2013the}, but solving them is computationally expensive; and
some use the Feynman-Kac formula to efficiently find the evolving expectation \cite{Chen2018Stochastic} but cannot obtain the probability distribution.
In addition, all these methods are \emph{intrusive} because they are not based on commercial simulation software. In these methods, the equations to be solved are derived through highly specialized mathematical knowledge, which have not been included in commercial software to date.

Seen from another perspective, as one of the state-of-the-art methods for uncertainty analysis that considers random parameters, polynomial chaos (PC) \cite{Xiu2010} has been introduced to power systems. In PC, a series expansion involving Wiener-Askey polynomials is used to represent the impact of random parameters on the system output. It precisely preserves the nonlinear nature of the system and is fast. Applications such as probabilistic power flow \cite{wu2017probabilistic,Ni2017Basis,sun2018probabilistic,xu2019probabilistic}, and dynamic uncertainty quantification \cite{Hockenberry2004,Li2017Uncertainty,prempraneerach2010uncertainty,Xu2019Propagating} have been proposed. Recently, due to community efforts, PC has become increasingly able to handle random parameters. For instance, adaptive sparse techniques have been developed to tackle high-dimensional input \cite{Ni2017Basis,sun2018probabilistic,Li2019Compressive,xu2019probabilistic}, and multi-element PC has been introduced to accommodate evolving probability densities in long-duration dynamic simulations \cite{prempraneerach2010uncertainty,Xu2019Propagating}.

However, although PC has shown its power in handling static random parameters, currently it has trouble dealing with continuous-time stochastic excitations \cite{debusschere2004numerical}. Karhunen-Lo\`{e}ve expansion (KLE) \cite{Park2015Modelling}, which spectrally decomposes a continuous random process into random parameters, may help PC handle stochastic excitations. Unfortunately, directly combining KLE and PC may result in tremendous complexity \cite{Park2015Modelling} and loss of information other than the mean and covariance in describing non-Gaussian uncertainty \cite{williams2006polynomial}. For this reason, existing work must assume that the excitations are Gaussian to facilitate the application of PC to their KLEs \cite{Li2019Compressive}.

To overcome the above discussed drawbacks of the existing methods to achieve decent dynamic uncertainty quantification of power systems subject to stochastic excitations, the following contributions are made in this paper: 1) the dynamic responses of power systems subject to non-Gaussian continuous disturbances are spectrally approximated as functions of discrete Gaussian random parameters based on SDE and KLE, and 2) by applying an adaptive sparse PC-based method to these functions, a method for quantifying the uncertainty in dynamic power systems subject to stochastic excitations that satisfies all three requirements raised above is then proposed.

This paper is organized as follows. Section \ref{sec:problem} introduces the problem of quantifying the power system dynamic uncertainty. Section \ref{sec:model} presents the It\^{o} process model of stochastic excitations, and the uncertainty quantification method is proposed in Section \ref{sec:method}. Finally, in Section \ref{sec:case}, case studies are used to verify the proposed method.

\section{Problem Description}
\label{sec:problem}

\begin{table}[tb]\scriptsize
  \renewcommand{\arraystretch}{1.35}
  \caption{Summary of the Features of Power System Dynamic Uncertainty Quantification Methods}
  \label{tab:methods}
  \centering
  \begin{tabular}{cccp{1.1cm}p{1.3cm}}
  \hline\hline
  Method        & Accuracy      & Efficiency      & Nonintrusive   & \hspace{2pt}\tabincell{c}{Continuous\\ Excitations}    \\
  \hline
  MCs \cite{dong2012numerical,milano2013systematic,ortega2018stochastic}          & High          & Low             & \hspace{12pt}Yes            &  \hspace{2pt}Applicable                     \\
  \tabincell{c}{SDE-Based   \cite{Apostolopoulou2016assessment,ju2018analytical,shi2018analytical,li2019analytic,Chen2018Stochastic}    }
                & Fair/Low      & High            & \hspace{13pt}No             & \hspace{2pt}Applicable                      \\
  \tabincell{c}{SDE-Based  \cite{Wang2013the}}
                & High          & Low             & \hspace{13pt}No             & \hspace{2pt}Applicable                   \\
  \tabincell{c}{PC-Based \cite{Hockenberry2004,prempraneerach2010uncertainty,Li2017Uncertainty,Xu2019Propagating}}
                & High          & High            & \hspace{12pt}Yes            & Inapplicable                        \\
  Proposed Method
                & High          & High           & \hspace{12pt}Yes            & \hspace{2pt}Applicable   \\
  \hline\hline
  \end{tabular}
\end{table}

For a power system subject to continuous-time disturbances, the dynamic response or performance index is a function of the entire \emph{time-domain realization} (or \emph{path}) of the stochastic excitations. Specifically, a power system subject to stochastic excitations can be modeled by a group of differential-algebraic equations with the stochastic excitations $\bm{\xi}_t$ as the parameter:
\begin{align}
  d{\bm{x}_t} &= \bm{f}(\bm{x}_t,\bm{y}_t;\bm{\xi}_t)dt, \label{eq:state} \\
  \bm{0} &= \bm{g} (\bm{x}_t,\bm{y}_t; \bm{\xi}_t),    \label{eq:algebra}
\end{align}
\noindent
where
$\bm{x}_t$ and $\bm{y}_t$ are the state and algebraic variables;
$t$ is time; and
$\bm{f}(\cdot)$ and $\bm{g}(\cdot)$ are the state and algebraic equations.

Without loss of generality, we assume that the solution of (\ref{eq:state})--(\ref{eq:algebra}) always exists. Then, by fixing the initial state $\bm{x}_0$, the system dynamic response, such as the trajectory of the rotor angle, or performance index, such as CPS1/2 in automatic generation control \cite{Chen2018Stochastic}, can be mathematically defined as a function of the entire path of the excitations $\bm{\xi}_t$, denoted by
\begin{align}
   \omega = \omega \big(\left\{\bm{\xi}_\tau\right\}_{\tau\in[0,t]} \big), \label{eq:response}
\end{align}
\noindent
where $\omega(\cdot)$ is referred to as the \emph{random response function (RRF)} in  this paper and $\omega$ represents the value of $\omega(\cdot)$.

The goal of quantifying the power system dynamic uncertainty is to extract statistical information on the RRF, such as the expectation and variance, which could help with evaluating the impact of renewable power integration and assessing the security in scheduling and operating the system.

If we replace the continuous excitation in (\ref{eq:state})--(\ref{eq:algebra}) by static random parameters, then the RRF (\ref{eq:response}) degrades to a function of the random parameters, and many mature methods, such as those based on polynomial chaos (PC), can be used to evaluate it. Unfortunately, when the excitations are continuous, no existing methods can simultaneously satisfy the three requirements, i.e., accuracy, efficiency, and nonintrusiveness, as noted in the Introduction and summarized in Table \ref{tab:methods}.

\begin{figure}[t]
  \centering
  \includegraphics[width=3.2in]{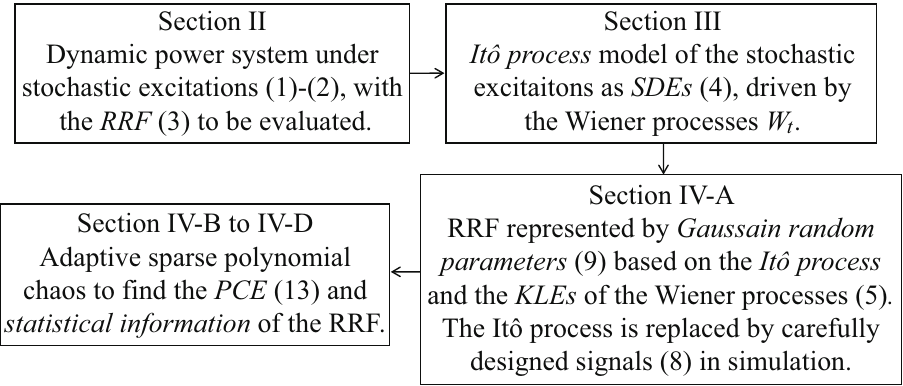}
  \caption{Brief framework of the proposed method.}
  \label{fig:frame}
\end{figure}

To fill this gap, this paper proposes a method for uncertainty quantification of power systems subject to stochastic excitations that satisfies all three of the requirements listed in the Introduction by combining the merits of methods based on SDEs and PC. Briefly, the Gaussian or non-Gaussian stochastic excitations in the power system are modeled with a It\^{o} process, and the RRF is represented by a 
function of static Gaussian random parameters. 
Then, the PCE of this function is calculated to extract statistical information about the RRF. The framework of the proposed method is shown in Fig. \ref{fig:frame}.

\section{Modeling the Stochastic Excitations as SDEs }
\label{sec:model}

\subsection{Modeling Stochastic Excitations Using the It\^{o} Process}
\label{sec:ito}

Many existing works have shown that the random processes of renewable generation can be represented by the It\^{o} process \cite{milano2013systematic,verdejo2016stochastic,Chen2018Stochastic}, formulated by the following SDE:
\begin{align}
  d\bm{\xi}_t & = \bm{\mu} \big( \bm{\xi}_t,t )dt + \bm{\sigma} \big( \bm{\xi}_t,t ) d{\bm{W}}_t,  \label{eq:ito}
\end{align}
\noindent
where $\bm{\xi}_t$ is the vector of stochastic excitations with dimension $m$; $\bm{W}_t = \left[ W_{1,t}, \ldots, W_{n,t} \right]^{\mathrm{T}}$ is the $n$-dimensional independent \emph{standard Wiener processes} (or \emph{Brownian motions}); and $\bm{\mu}(\cdot,\cdot): \mathbb{R}^m \times \mathbb{R}_+ \rightarrow \mathbb{R}^m$ and $\bm{\sigma}(\cdot,\cdot): \mathbb{R}^m \times \mathbb{R}_+ \rightarrow \mathbb{R}^{m \times n}$ are called the \emph{drift} and \emph{diffusion terms}, respectively.

\begin{remark}
\label{remark:1}
With different drift and diffusion terms, Gaussian or non-Gaussian continuous random processes $\bm{\xi}_t$ with arbitrary probability distributions can be depicted by (\ref{eq:ito}) \cite{Chen2018Stochastic}. Therefore, modeling continuous random disturbances in power systems would not be a limitation of the proposed method.
\end{remark}

For example, when the drift term $\bm{\mu}(\cdot,\cdot)$ is affine to $\bm{\xi}_t$ and the diffusion term $\bm{\sigma}(\cdot,\cdot)$ is constant, $\bm{\xi}_t$ becomes the \emph{Ornstein-Uhlenbeck process}, which is Gaussian \cite{Pardoux2014Stochastic} and used in many other works \cite{Apostolopoulou2016assessment,ju2018analytical,shi2018analytical,li2019analytic}. Also, using the drift and diffusion terms in Table \ref{tab:ito}, typical Non-Gaussian excitation can be modeled.

Noted that in Table \ref{tab:ito} the drift and diffusion terms do not have an explicit time argument. This is because the listed It\^{o} processes are time-homogeneous, and hence the drift and diffusion terms are not explicit functions of time. In contrast, in the reminder of this paper, a time argument explicitly appears in the drift and diffusion terms for generality. The proposed uncertainty quantification method is able to handle both homogeneous and inhomogeneous stochastic excitations.

The relation between the probability distribution and the drift and diffusion terms is determined by the \emph{Fokker-Planck equation} \cite{Pardoux2014Stochastic,Wang2013the} (see Appendix \ref{sec:fp}), and the It\^{o} process can be analytically constructed based on the it; see \cite{Chen2018Stochastic} for details. Recent works also provide method to constructing the It\^{o} process model for renewable energies based on analytical model \cite{zarate2016construction} or real data \cite{jonsdottir2019data} with arbitrary probability distribution and temporal correlation.

In addition, we present a data-driven method for identifying the It\^{o} process model from real data; see Appendix \ref{sec:data}.

Although modeling continuous random excitations is not a limitation of this work as noted in Remark \ref{remark:1}, it needs to be clarified that random jumping disturbances in power systems, e.g., random load switchings, cannot be precisely modeled by the It\^{o} process, as the It\^{o} process is driven by Wiener processes that are continuous. The Poisson process or L\'{e}vy process could be more suitable for describing random jumping disturbances \cite{dong2012numerical,mele2019modeling}. Because this paper focuses on analyzing power systems under continuous disturbances, jumping disturbances will not be discussed further.

\begin{table}[tb]\scriptsize
  \renewcommand{\arraystretch}{1.6}
  \caption{Drift and Diffusion Terms of One-Dimensional It\^{o} Processes with Some Typical Probability Distributions}
  \label{tab:ito}
  \centering
  \begin{tabular}{cccc}
  \hline\hline
  Type          & Probability Density                                               & $\mu\left(\xi_t\right)$                & $\sigma^2\left(\xi_t\right)$        \\
  \hline
  Gaussian      & \tabincell{c}{$e^{{(\xi_t-a)^2}/2b}/\sqrt{2 \pi b}$}              & $-\left(\xi_t-a\right)$                & $2b$                  \\
  Beta          & \tabincell{c}{$\dfrac{\xi_t^{a-1}(1-\xi_t)^{b-1}}{B(a,b)}$}       & $-\left(\xi_t-\dfrac{a}{a+b}\right)$   & $\dfrac{2\xi_t\left(1-\xi_t\right)}{a+b}$        \\
  Gamma         & \tabincell{c}{$\dfrac{b^a}{\Gamma(a)}\xi_t^{a-1}e^{-b\xi_t}$}     & $-\left(\xi_t-a/b\right)$              & $2\xi_t/b$            \\
  Laplace       & \tabincell{c}{$e^{-|\xi_t-a|/b}/(2b)$}                            & $-\left(\xi_t-a\right)$                & $2b|\xi_t-a|+2b^2$    \\
  \hline\hline
  \end{tabular}
\end{table}

\subsection{Simulation of the Stochastic Excitations}

Monte Carlo simulation (MCs) is the most straightforward way to quantify the statistical information in the RRF (\ref{eq:response}) \cite{dong2012numerical,milano2013systematic,ortega2018stochastic}.
In MCs, paths of the stochastic excitations are sampled via a stochastic numerical integration scheme, for example the Maryuama-Euler (EM) scheme \cite{milano2013systematic}, as
\begin{align}
   \bm{\xi}_{t + h} & = h  \bm{\mu} \left(  \bm{\xi}_t, t \right) +  \bm{\sigma} \big( \bm{\xi}_t, t ) \sqrt{h}\bm{\zeta}, \label{eq:distito}
\end{align}
\noindent
where $h$ is the step length; $\bm{\zeta} \sim \mathcal{N}\left(\bm{0},\bm{I} \right)$ is sampled as a vector of independent normal random variables in each step.

Integration scheme (\ref{eq:distito}) implies time-domain discretization of the stochastic excitations is impractical in uncertainty quantification, since the number of the resulting random parameters is proportional to the number of discretization steps, causing the curse of dimensionality. However, discretizing the excitations spectrally instead of in the time domain can be a feasible solution. This is the basic idea of our work, as illustrated later.

\section{The Uncertainty Quantification Method}
\label{sec:method}

\subsection{Spectral Representation of the Stochastic Excitations}
\label{sec:kl}

According to the Karhunen-Lo\`{e}ve theorem \cite{friz2010multidimensional}, a standard Wiener process $W_{i,t}$ can be orthogonally decomposed into a series of independent standard Gaussian random variables $\{\zeta_{i,j}\}_{j=1}^{\infty}$, known as the \emph{Karhunen-Lo\`{e}ve expansion (KLE)}:
\begin{align}
  W_{i,T} = \int_0^T dW_{i,t} = \sum_{j=1}^{\infty} \zeta_{i,j} \int_0^T m_j(t) dt, \label{eq:kl}
\end{align}
\noindent
where $\left\{m_j(t)\right\}_{j=1}^{\infty}$ are functions defined on interval $t\in[0,T]$:
\begin{align}
{m_j(t)=}
\begin{cases}
    \sqrt{{1}/{T}} &,\ j = 1,\\
    \sqrt{{2}/{T}} \cos \left[{(j-1)\pi t}/{T} \right] &,\ j \ge 2. \label{eq:klb2}
\end{cases}
\end{align}

For computation, the infinite series (\ref{eq:kl}) is truncated at a given order $K$. 
Taking the derivative of both sides of (\ref{eq:kl}) yields
\begin{align}
    dW_{i,t} \approx \sum_{j=1}^{K}  \zeta_{i,j} m_j(t). \label{eq:kle}
\end{align}

Substituting (\ref{eq:kle}) into the It\^{o} process (\ref{eq:ito}) yields the following ordinary differential equation with random coefficients:
\begin{align}
     \frac{d\bm{\xi}^*_t}{dt} (\bm{\zeta})& = \bm{\mu} \big( \bm{\xi}^*_t,t \big) + \bm{\sigma} \big( \bm{\xi}^*_t,t \big)  \sum_{j=1}^{K}  \bm{\zeta}_{j} m_j(t), \label{eq:eqkl}
\end{align}
\noindent
where $\bm{\zeta}_{j}$ is a vector of independent normal random variables with dimension $n$ and $\bm{\zeta}\triangleq \{\zeta_{i,j}\}_{1 \le i \le n,1 \le j \le K}$ is the vector of all independent normal random variables. For convenience, in the rest of this paper, we reindex the entries of $\bm{\zeta}$ as $\{\zeta_i\}_{1\le i \le M}$ without ambiguity. $M=nK$ is the size of $\bm{\zeta}$.

The truncation order $K$ can be empirically chosen from $3$ to $6$ based on a compromise between precision and computing cost.
Although currently we have not yet found the analytical relation between $K$ and the accuracy of the proposed method, intuitively, we can see that $K$ roughly determines the frequencies of the alternating components of the excitation signals (\ref{eq:eqkl}). When evaluating the uncertainty in power system electromechanical dynamics over an interval of about ten seconds, such choice makes the frequencies of the alternating components of the excitation signal (\ref{eq:eqkl}) match the time constant of power system dynamics, ensuring the proposed method to obtain main information. Numerical results in Sections \ref{sec:39} and \ref{sec:118} also show that such choice is acceptable.

Visual examples of $\bm{\xi}^*_t(\bm{\zeta})$ as the solution of (\ref{eq:eqkl}) for fixed values of $\bm{\zeta}$ are shown in Fig. \ref{fig:cps}(b). An intuitive proposition regarding (\ref{eq:eqkl}) is as follows:
\begin{proposition}
\label{prop:1}
The solution of (\ref{eq:eqkl}), i.e., $\bm{\xi}^*_t(\bm{\zeta})$, converges to the It\^{o} process $\bm{\xi}_t$ defined in (\ref{eq:ito}) as $K \rightarrow \infty$.
\end{proposition}
Although this proposition seems straightforward, the rigorous proof is rather arduous. Interested readers are referred to Section 15.5.3 of \cite{friz2010multidimensional}. With such an approximation, the non-Gaussian stochastic excitations can be characterized by independent Gaussian random variables, which directly facilitates the application of the polynomial chaos based methods.

Finally, by substituting (\ref{eq:eqkl}) into (\ref{eq:response}), we find that the system dynamic response or performance index defined by the RRF $\omega(\cdot)$ can be approximated by an implicit function (denoted by $\omega^*(\cdot)$) of the independent Gaussian random variables, as
\begin{align}
    \omega \approx \omega^*(\bm{\zeta}) = \omega\big(\left\{\bm{\xi}^*_\tau(\bm{\zeta})\right\}_{\tau\in[0,t]} \big). \label{eq:implict}
\end{align}

Analysis of the convergence of approximation (\ref{eq:implict}) is delayed to the proof of Proposition \ref{prop:2}.

This far, we have established the foundation for using Gaussian random parameters to depict the RRF during stochastic excitation as (\ref{eq:implict}), but $\omega^*(\bm{\zeta})$  has still not been expressed explicitly. Fortunately, polynomial chaos provides a way to approximate it explicitly. This is illustrated followingly.

\subsection{Quantifying Uncertainty Using Polynomial Chaos}
\label{sec:pc}

Polynomial chaos (PC) uses a series of polynomials orthogonal with respect to the probability density of the random parameters to depict the random outputs \cite{Xiu2010}. For the standard Gaussian random variables in the KLEs (\ref{eq:kl}) of the Wiener process, the relevant polynomials are the Hermite polynomials,
\begin{align}
    H_n(\zeta_{i}) = (-1)^n e^{\zeta_{i}^2/2} \frac{d^n}{d\zeta_{i}^n} e^{-\zeta_{i}^2/2},\ n\in\mathbb{N},
\end{align}
\noindent
which admit the following orthogonality:
\begin{align}
    \left\langle H_j(\cdot),  H_k(\cdot) \right\rangle & \triangleq \int_{-\infty}^{+\infty} H_j  H_k  \frac{e^{-\zeta_i^2/2}}{\sqrt{2\pi}}d\zeta_i
    = \delta_{jk} \| H_j \|^2,
\end{align}
\noindent
where $\langle\cdot,\cdot\rangle$ represents an inner product; $\|\cdot\|=\langle\cdot,\cdot\rangle^{1/2}$ is the $2$-norm; and $\delta_{jk}$ is the Kronecker delta function.

Then, given an order $N_i$ for each scalar random variable ${\zeta}_i$, the basis $\left\{ \phi_{\bm{j}}(\bm{\zeta}) \right\}$ for the random vector $\bm{\zeta}$ is constructed as
\begin{align}
    \hspace{-2pt}\left\{ \phi_{\bm{j}}(\bm{\zeta}) \right\} = \left\{ H_{j_1}(\zeta_1) H_{j_2}(\zeta_2) \cdots H_{j_M}(\zeta_{j_{M}}) : j_i \le N_i \right\}, \label{eq:wick}
\end{align}
\noindent
where $\bm{j} \triangleq \left[ j_1,j_2,\ldots,j_{M} \right]^{\mathrm{T}}$ is the multi-dimensional index.



Finally, the approximate RRF $\omega^*(\bm{\zeta})$ in (\ref{eq:implict}) can be explicitly approximated by the following polynomial chaos expansion (PCE), which is denoted $\hat{\omega}(\cdot)$ and formulated as follows:
\begin{align}
    \omega \approx \omega^*(\bm{\zeta}) \approx \hat{\omega}(\bm{\zeta}) \equiv \sum_{\bm{j}} \hat{c}_{\bm{j}} \phi_{\bm{j}}(\bm{\zeta}), \label{eq:pce}
\end{align}
\noindent
where $\hat{c}_{\bm{j}}$ is the coefficient of the PCE. Its exact value is defined by the orthogonal projection $\hat{c}_{\bm{j}} = \langle \omega^*(\bm{\zeta}), \phi_{\bm{j}}(\bm{\zeta}) \rangle / {\|\phi_{\bm{j}}(\bm{\zeta})\|^2}$ \cite{Xiu2010} but is generally unavailable due to the absence of an explicit expression for $\omega^*(\bm{\zeta})$.

Instead of rigid orthogonal projection, in practice, either the probabilistic Galerkin method (PGM) or the probabilistic collocation method (PCM) can be used to find the coefficients \cite{Xiu2010}. Because the PCM is a nonintrusive method whereas the PGM is not, this study uses the PCM to find the coefficients.

\subsection{PCM for Computing the PCE Coefficients}
\label{sec:pcm}

First, define the set of collocation points (Gaussian quadrature points) that contain all the zeros of the product of the $(N_i+1)$th Hermite polynomials of the random variables in $\bm{\zeta}$:
\begin{align}
    \big\{ \hat{\bm{\zeta}}_i : H_{N_1+1}(\zeta_1) H_{N_2+1}(\zeta_2) \cdots H_{N_{M}+1}(\zeta_M) = 0 \big\}_{i=1}^{N_b}, \label{eq:cps}
\end{align}
\noindent
where $N_b$ is the number of collocation points, which exactly equals the number of basis functions in (\ref{eq:wick}).

For each collocation point $\hat{\bm{\zeta}}_i$ in (\ref{eq:cps}), (\ref{eq:eqkl}) is used to generate the path $\bm{\xi}^*_t(\hat{\bm{\zeta}}_i)$ of excitation $\bm{\xi}_t$ as the external input signal, and dynamic simulation software 
is used to compute the corresponding response,
$\hat{\omega}_i=\omega \big( \{ \bm{\xi}^*_\tau(\hat{\bm{\zeta}}_i) \}_{\tau \in [0,t]}  \big)$. Then, the PCE coefficients in (\ref{eq:pce}) can be computed using
\begin{align}
    \left[  \hat{c}_1, \hat{c}_2 , \ldots, \hat{c}_{N_b}  \right]^{\mathrm{T}} =
    \bm{A}^{-1}
    \left[  \hat{\omega}_1, \hat{\omega}_2 , \ldots, \hat{\omega}_{N_b}  \right]^{\mathrm{T}}, \label{eq:cm}
\end{align}
\noindent
where $\bm{A}$ is a constant matrix  solely determined by the basis (\ref{eq:wick}) (reindexed as $\left\{ \phi_i(\cdot) \right\}_{i=1}^{N_b}$) and the collocation points (\ref{eq:cps}):
\begin{align}
    \bm{A} =
    \begin{bmatrix}
        \phi_1(\hat{\bm{\zeta}}_1)     & \cdots & \phi_{N_b}(\hat{\bm{\zeta}}_1) \\
        \vdots                       & \ddots & \vdots   \\
        \phi_1(\hat{\bm{\zeta}}_{N_b}) & \cdots &  \phi_{N_b}(\hat{\bm{\zeta}}_{N_b})
         \\
    \end{bmatrix}.  \label{eq:cmat}
\end{align}

Further, since only a few higher-order terms have a noticeable impact on the precision of the PCE \cite{Xiu2010}, we do not need to include all the higher-order terms to construct the PCE.

Alternatively, we use the \emph{Smolyak adaptive sparse algorithm}\cite{CONRAD2013ADAPTIVE} to incrementally add higher-order terms in an adaptive way until a given error tolerance $\epsilon$ is reached. The precision of the PCE is assessed by the coefficients. Only the terms with a noticeable impact are added, which markedly decreases the computational cost.

The Smolyak algorithm ensures the error in the $2$-norm (variance in this case) of the calculated PCE generally has at most a same order of magnitude as $\epsilon$. This means if we need a $1\%$-level accuracy in evaluating the statistical information on the RRF, choosing $\epsilon$ to be $0.01$ or $0.001$ is sufficient. Since this algorithm is not a contribution of this paper, to save space, we do not discuss it further. Instead, details of the procedure can be found in \cite{CONRAD2013ADAPTIVE}.

\subsection{Mathematical Foundation of the Proposed Method}
\label{sec:conv}

Before laying the mathematical foundation of the proposed method, we make several generic assumptions.

\begin{assumption}
\label{assu:2}
During the dynamic process, the number of triggering conditions (denoted as $s_i(\bm{x}_t,\bm{y}_t;\bm{\xi}_t)=0,\ i\in\Omega^\mathrm{S}$, where $\Omega^\mathrm{S}$ is the set of switching conditions; see \cite{Hiskens2000}) of switching events, e.g., of on-load tap changers and over-excitation limiters, is finite. Because in practice power systems are analyzed over a finite time interval $\tau\in[0,t]$, this assumption naturally holds.
\end{assumption}

\begin{assumption}
\label{assu:3}
The RRF is finite, i.e. $| \omega\left( \cdot \right) | \le \infty$, almost surely over the finite time interval $\tau\in[0,t]$. This assumption is guaranteed as long as no singularity occurs in the power system differential-algebraic model (see \cite{wang2011investigation}).
\end{assumption}

\begin{assumption}
\label{assu:1}
The RRF has the Lipschitz continuity in terms of $L^2_w$-norm, i.e., there exists a constant $L$ such that for arbitrary two stochastic excitations, say $\{\bm{\xi}^{(1)}_\tau\}_{\tau\in[0,t]}$ and $\{\bm{\xi}^{(2)}_\tau\}_{\tau\in[0,t]}$, we have $\big\| \omega\big(\{\bm{\xi}^{(1)}_\tau\}_{\tau\in[0,t]}\big) - \omega\big(\{\bm{\xi}^{(2)}_\tau\}_{\tau\in[0,t]}\big) \big\| \le L \big\| \{\bm{\xi}^{(1)}_\tau - \bm{\xi}^{(2)}_\tau \}_{\tau\in[0,t]} \big\|$, where $\|\cdot\|$ is the $L^2_w$-norm of the square Lebesgue integrable space $L^2_w$ equipped with probability measure $w$.
\end{assumption}

Then, we have the following proposition laying the mathematical foundation of the proposed method:

\begin{proposition}
\label{prop:2}
The PCE (\ref{eq:pce}) converges to the RRF (\ref{eq:response}), i.e., $\| \omega\left( \{ \bm{\xi}_\tau \}_{\tau\in[0,t]}\right) - \hat{\omega}\left( \bm{\zeta} \right) \| \rightarrow 0$ , as $K \rightarrow \infty$ and the error tolerance of the Smolyak adaptive sparse collocation method $\epsilon \rightarrow 0$. See proof in Appendix \ref{sec:proof}.
\end{proposition}

Note that the existence of discrete switching events do not affect Assumption \ref{assu:1} and then Proposition \ref{prop:2}. This is, not rigidly speaking, because switching events only impact the $L^2_w$-norm when the triggering manifold $\cup_{i\in\Omega^\mathrm{S}} \big(s_i(\bm{x}_t,\bm{y}_t;\bm{\xi}^*_t(\bm{\zeta}))=0\big)$ have a non-zero measure in the probability space of $\bm{\zeta}$, but Assumption \ref{assu:2} rules out this situation.

However, extreme nonlinearity or discrete switching events in power system dynamics may negatively impact the convergence rate of the proposed method, because they deteriorate the smoothness of the RRF, as quantified in (\ref{eq:appconv}). This deterioration may lead to inaccurate results in certain situations and currently is a limitation of the proposed method. An extreme case is shown in Section \ref{sec:instability}.

Note that switching events do not necessarily cause such deterioration. In Section \ref{sec:case} we impose and clear grounding faults in all cases, which are discrete switching events. Because they do not cause discontinuity in the RRF, the proposed method achieves accurate results. The only case causing discontinuity in the RRF is that with different paths of the stochastic excitations, some system dynamic trajectories encounters a switching condition and then jumpings in occurs the system states, whereas other trajectories remain continuous. This causes the trajectories become two separate clusters, which is visually similar to Fig. \ref{fig:traj}, reflecting the RRF is no more continuous with respect to the stochastic excitations.

To alleviate such limitation, multi-element PC may be a feasible solution. In multi-element PC, the support domains of random variables are divided into multiple elements, within each element the output is relatively smooth with respect to the inputs, and therefore the convergence rate is preserved \cite{prempraneerach2010uncertainty,Xu2019Propagating}. Leveraging multi-element PC in the framework of the proposed method is worth investigating in the future study.

\begin{figure}[tb]
  \centering
  \includegraphics[width=3.5in]{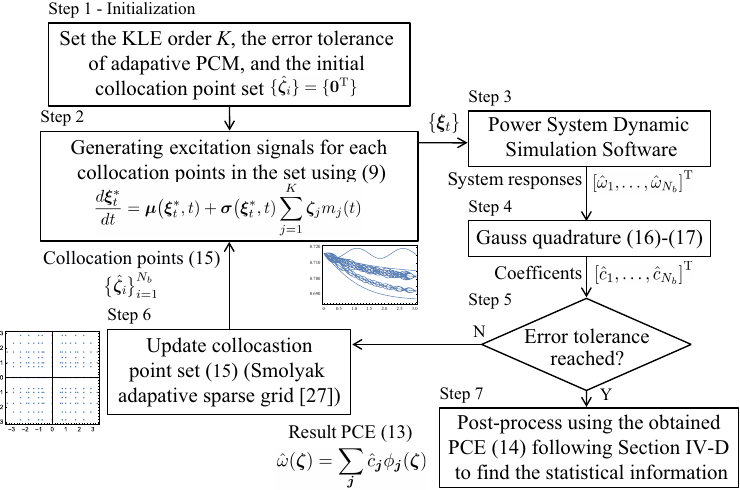}
  \caption{Flow chart of the proposed uncertainty quantification method.}
  \label{fig:flow}
\end{figure}

\subsection{Post-Processing to Extract Statistical Information}
\label{sec:post}

From the coefficients in the PCE (\ref{eq:pce}), we can easily find the expectation and variance of the RRF \cite{Xiu2010} as follows:
\begin{align}
  \mathbb{E} & \left[ \omega(\cdot) \right] = \hat{c}_{\bm{0}},\  \mathbb{E} \left[ \omega^2(\cdot) \right] = \sum_{\bm{j}} \hat{c}^2_{\bm{j}} \gamma_{\bm{j}},
\end{align}
\noindent
where $\hat{c}_{\bm{0}}$ is the coefficient of the constant basis function in (\ref{eq:wick}) with all $j_i=0$, and $\gamma_{\bm{j}} = \|\phi_{\bm{j}}(\cdot)\|^2 = j_1!j_2!\cdots j_n!$.

Additionally, MCs can be used with the PCE (\ref{eq:pce}) to obtain statistical information, e.g., high-order moments. This can be very efficient since the PCE is a simple explicit function and does not require time-domain simulation of the power system.

\subsection{Overall Computational Procedure}

\begin{algorithm}[tb]
  \caption{Procedure of the Proposed Uncertainty Quantification Method}
  \label{alg:pipm}
  \begin{algorithmic}[1]
    \REQUIRE It\^{o} process model (\ref{eq:ito}) of the stochastic excitations, truncation order $K$ of the KLE (\ref{eq:kle}), error tolerance $\epsilon$ of the Smolyak collocation method, dynamic power system model in the simulation software
    \STATE \label{step:1} initialize the set of collocation points $\{ \hat{\bm{\zeta}}_i \}$ with only one element $\bm{0}^{\mathrm{T}}$
    \STATE \label{step:2} for each sampling point $\hat{\bm{\zeta}}_k$, use (\ref{eq:eqkl}) to calculate the paths of corresponding disturbances
    \STATE \label{step:3} use power system simulation software, such as PSS/E, to compute the values of the RRF for each external excitation signal obtained in Step 2
    \STATE \label{step:4} calculate the PCE coefficients based on the PCM, specifically the Gaussian quadrature (\ref{eq:cm})--(\ref{eq:cmat})
    \STATE \label{step:5} based on the change in the coefficients, estimate approximation error using the method in \cite{CONRAD2013ADAPTIVE}; if error tolerance $\epsilon$ is reached, go to Step 7; otherwise, go to Step 6
    \STATE \label{step:6} update the set of collocation points $\{ \hat{\bm{\zeta}}_i \}$ based on the changes in the coefficients \cite{CONRAD2013ADAPTIVE}, and return to Step 2
    \STATE \label{step:7} calculate the statistics using the resulting PCE following Section \ref{sec:post}
  \end{algorithmic}
\end{algorithm}

The overall procedure of the proposed method is summarized as Algorithm \ref{alg:pipm}, and the flow chart is shown as Fig. \ref{fig:flow}.

\section{Case Studies}
\label{sec:case}

\subsection{It\^{o} Process Modeling of Empirical Stochastic Excitations and Analytical Models}
\label{sec:casemodel}

We start by demonstrating the ability of the It\^{o} process model (\ref{eq:ito}) to model the uncertainty in renewable, for example wind and solar power generation on different time scales, based on empirical data and an analytical uncertainty model.

\subsubsection{Modeling Wind Power Uncertainty over Minutes Based on Empirical Data}

\begin{figure}[t]
  \centering
  \includegraphics[width=2.580in]{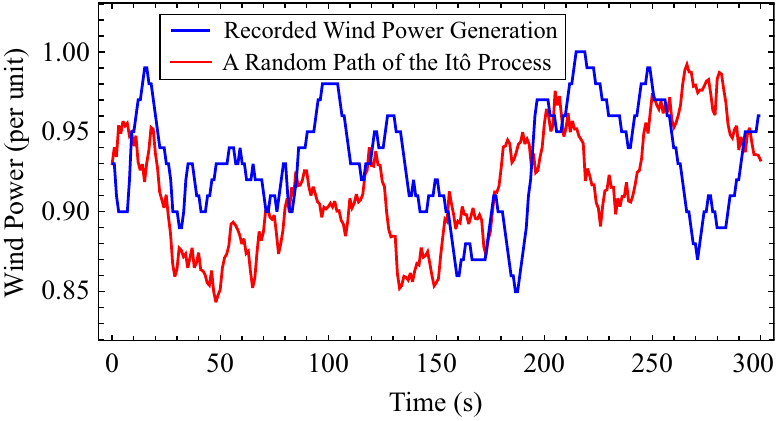}
  \caption{Power generation of an offshore wind farm over a $5$-minute interval recorded by Ris\o\, and a random path of the identified It\^{o} process model.}
  \label{fig:wind} \vspace{12pt}
  \includegraphics[width=3.4in]{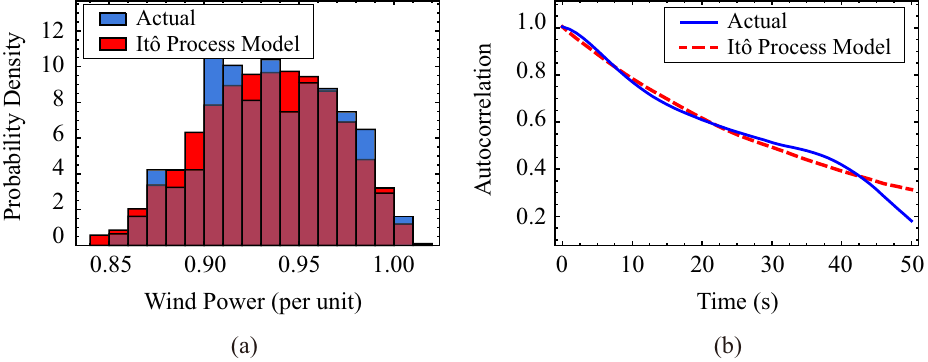}
  \caption{Probability density and autocorrelation of the actual wind power and the identified It\^{o} process model. (a) Probability density. (b) Autocorrelation.}
  \label{fig:wid}
\end{figure}

Fig. \ref{fig:wind} shows the per-unit electrical power generated by an offshore wind farm over a $5$-minute interval recorded by Ris\o\ DTU at a resolution of $1$ point per second \cite{lin2012assessment}. The probability density is shown in Fig. \ref{fig:wid}(a), and the autocorrelation is shown in Fig. \ref{fig:wid}(b).

As can be seen, the distribution density of the recorded data shows a complicated shape. Although based on the Fokker-Planck equation (\ref{eq:fp}) we can construct the It\^{o} process with exactly the same probability density \cite{Chen2018Stochastic}, it leads to a weird-shaped diffusion term indicating overfitting and can be inappropriate for practical application. To avoid overfitting, we use quadratic polynomials to construct the drift and diffusion terms. Using the method introduced in the Appendix, the identified It\^{o} process model is formulated as follows:
\begin{align}
  dP_t = & \left[ 0.0535 - 0.0899 P_t + 0.0349 P_t^2 \right] dt  \nonumber \\
  & +  \left[-0.410 + 0.919 P_t - 0.505 P_t^2 \right] dW_t. \label{eq:windito}
\end{align}

A random path of the identified model (\ref{eq:windito}) is simulated and shown in Fig. \ref{fig:wind}. Visually, the path has a shape of volatility similar to that of the actual curve. Quantitative comparisons of the probability distribution and the autocorrelation are shown in Fig. \ref{fig:wid}. Clearly, the non-Gaussian distribution and temporal correlation are well characterized.

\subsubsection{Modeling Solar Power Uncertainty Over One Day Based on Empirical Data}

\begin{figure}[t]
  \centering
  \includegraphics[width=2.58in]{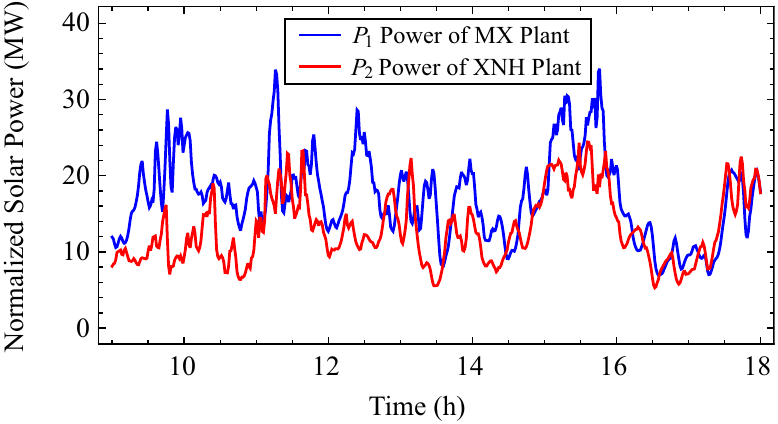}
  \caption{Solar power output on May 30, 2018 at two pv plants in Sichuan, China. Both are normalized by the sine of the solar altitude angle.}
  \label{fig:solar} \vspace{12pt}
  \includegraphics[width=3.3in]{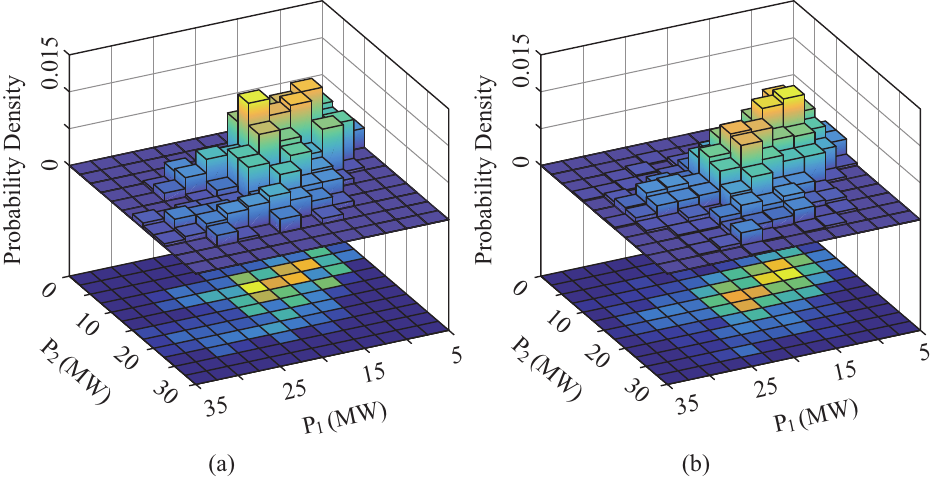}
  \caption{Joint probability densities of (a) the recorded solar power output, and (b) the simulation of the corresponding It\^{o} process model. }
  \label{fig:identify}
\end{figure}

Fig. \ref{fig:solar} shows the electrical power generated and injected to the grid by two nearby photovoltaic plants in Sichuan, China on May 30, 2018. Both are normalized by the sine of the solar altitude angle. Their joint probability density is shown in Fig. \ref{fig:identify}(a).

Based on Fig. \ref{fig:identify}(a), we assume that there are two correlated beta distributions. Then, the It\^{o} process model is identified as
\begin{align}
  d\left[
    \begin{matrix}
            P_1    \\
            P_2 \\
    \end{matrix} \right]
  = &
    \left[
    \begin{matrix}
            0.00500 - 0.000280 P_1    \\
            0.00513 - 0.000374 P_2 \\
    \end{matrix} \right] dt \nonumber \\& +
      \left[
    \begin{matrix}
            0.00442 \sqrt{a} & 0  \\
            0.00271 \sqrt{a} & 0.00364 \sqrt{b} \\
    \end{matrix} \right]
  \left[
    \begin{matrix}
            dW_{1,t}  \\
            dW_{2,t} \\
    \end{matrix} \right],
\end{align}
\noindent
where $a = -598.5 + 124.7 P_1 - P_1^2$; $b = -404.3 + 85.9 P_2 - P_2^2$. The non-diagonal entries in the diffusion term represent the correlation between the power generated by the two plants.

The simulated joint probability density of the It\^{o} process model is shown in Fig. \ref{fig:identify}(b). In comparison with Fig. \ref{fig:identify}(a), the model well depicts the probability distribution and reflects the correlation between the two plants.

\subsubsection{Modeling Wind Speed Uncertainty Based on the Analytical Model}
\label{sec:analytical}

In engineering practices, the utilities often use typical analytical uncertainty models to represent the uncertainty of renewable generations. The community has shown that these analytical model can be exactly represented by It\^{o} processes \cite{zarate2016construction,jonsdottir2019data}. Based on the analytical model, we show that the method presented in Appendix \ref{sec:data} can accurately construct the It\^{o} process model based on simulation data.

For wind speed (denoted as $v_t$) that admits a two-parameter Weibull distribution and an exponential temporal correlation, the probability density function $p(v_t)$ and autocorrelation function $\rho_{v_t}(\tau)$  \cite{zarate2016construction} are
\begin{align}
 p(v_t) &=
    {\lambda_1}/{\lambda_2}\left( {v_t}/{\lambda_2} \right)^{\lambda_1-1} \times \exp \big[ -\left( {v_t}/{\lambda_2} \right)^{\lambda_1} \big], \label{eq:windpdf} \\
  \rho_{v_t}(\tau) &= \exp[-\alpha\tau] \label{eq:windautoco}
\end{align}
\noindent
where $\lambda_1$, $\lambda_2$, and $\alpha$ are parameters. In this case, we set $\lambda_1=2.343$, $\lambda_2=5.244$, and $\alpha=0.25$ \cite{zarate2016construction}.

According to \cite{zarate2016construction}, the It\^{o} process can be exactly characterized using the following drift and diffusion terms, as
\begin{align}
  \mu(v_t) =& -\alpha \left[ v_t - \lambda_2 \Gamma \left(1+1/{\lambda_1} \right) \right] \label{eq:weidrift} \\
  \sigma(v_t) =&  \sqrt{\sigma_1(v_t)\sigma_2(v_t)} \label{eq:weidiff} \\
    \text{with}\hspace{40pt}& \hspace{380pt} \nonumber \\
  \sigma_1(v_t)=& 2\alpha ({\lambda_2}/{\lambda_1^2}) v_t \left( {\lambda_2}/{v_t} \right)^{\lambda_1}, \nonumber \\
  \sigma_2(v_t)=& \lambda_1 \exp\big[ \left( {v_t}/{\lambda_2}\right)^{\lambda_1} \big] \Gamma\big( 1+{1}/{\lambda_1},\left( {v_t}/{\lambda_2} \right)^{\lambda_1}\big)  \nonumber \\ & - \Gamma\big({1}/{\lambda_1}\big). \nonumber
\end{align}
\noindent
where $\Gamma(\cdot)$ is the gamma function.

\begin{figure}[t]
  \centering
  \includegraphics[width=3.48in]{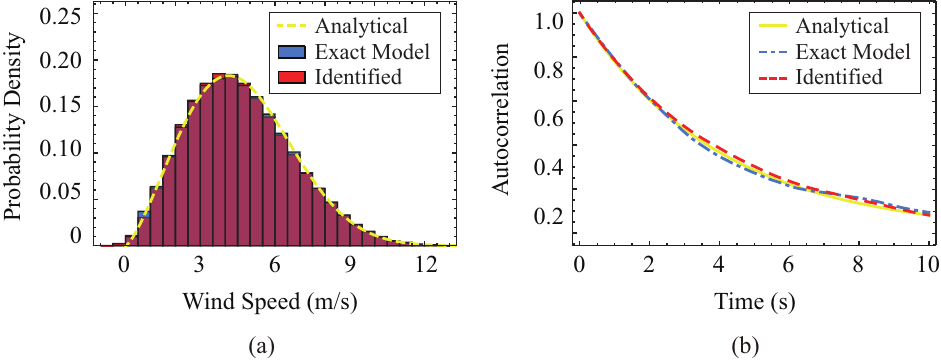}
  \caption{The probability density and autocorrelation functions of the exact and identified It\^{o} process models, compared to the analytical curves. (a) Probability density function. (b) Autocorrelation function.}
  \label{fig:modelcheck}
\end{figure}

In addition, based on sampled paths of the exact model (\ref{eq:weidrift})-(\ref{eq:weidiff}), assuming linear drift term and quadratic diffusion term, the following It\^{o} process model is identified:
\begin{align}
  dv_t = & \left[ 1.122 -0.243 v_t  \right] dt  \nonumber \\
  & +  \left[0.763 + 0.224 v_t - 0.016 v_t^2 \right] dW_t. \label{eq:windito}
\end{align}

By the MCs, the density and autocorrelation functions of the identified It\^{o} process model (\ref{eq:windito}) are compared to the exact It\^{o} model (\ref{eq:weidrift})--(\ref{eq:weidiff}) and the analytical model (\ref{eq:windpdf})--(\ref{eq:windautoco}) in Fig. \ref{fig:modelcheck}.

As the result, the density and autocorrelation functions almost overlap with each other, showing that the same as the analytical model, the identified model is also capable of precisely depicting the stochastic excitations.

\subsection{IEEE 39-Bus System Case Study}
\label{sec:39}

\subsubsection{Setting}

The proposed uncertainty quantification method is first tested on the IEEE 39-bus system \cite{39sys}.
The method is implemented in \emph{Python} with dynamic simulation performed by \emph{PSS/E} via \emph{PSSPY} interface. Detailed models of GENROU generators, IEEET1 exciters, and TGOV1 governors are included. A desktop with an Intel i7-8700 CPU is used.

To highlight the ability of the proposed method to handle different types of non-Gaussian uncertainty, three power injections $P_3, P_{15}$ and $P_{29}$ with correlated beta and Laplace distributions are connected to buses 3, 15, and 29 to model volatile power generation. Their per-unit values with a base power of $100$ MW are modeled by the following It\^{o} process:
\begin{align}
  d\left[
  \setlength{\arraycolsep}{4pt}
  \renewcommand{\arraystretch}{0.95}
  \begin{matrix}
    P_3 \\ P_{15} \\ P_{29}
  \end{matrix} \right]
  = &
  \left[
  \setlength{\arraycolsep}{4pt}
  \renewcommand{\arraystretch}{0.95}
  \begin{matrix}
    0.268 - 0.080 P_3    \\
    0.160 - 0.050 P_{15} \\
    0.270 - 0.100 P_{29}
  \end{matrix} \right]
  dt  \nonumber \\
  & +
  \left[
  \setlength{\arraycolsep}{4pt}
  \renewcommand{\arraystretch}{0.95}
  \begin{matrix}
    0.163 \sqrt{c}  & 0.1  & 0 \\
    -0.082 \sqrt{c} & 0.2  & 0 \\
    0               & 0.1  & 0.2 \sqrt{d}
  \end{matrix} \right]
  \left[
  \setlength{\arraycolsep}{4pt}
  \renewcommand{\arraystretch}{0.95}
  \begin{matrix}
    dW_{1,t} \\ dW_{2,t} \\ dW_{3,t}
  \end{matrix} \right], \label{eq:case}
\end{align}
\noindent
where $c = -10 + 6.5 P_3 - P_3^2$; $d= |2.73 - P_{29}| + 0.01$. To further demonstrate nonlinearity, a three-phase grounding fault is imposed on bus 3 at $t=1$ s and lasts for $0.2$ s until it is cleared by tripping line 3-4.

\begin{figure}[tb]
  \centering
  \includegraphics[width=3.42in]{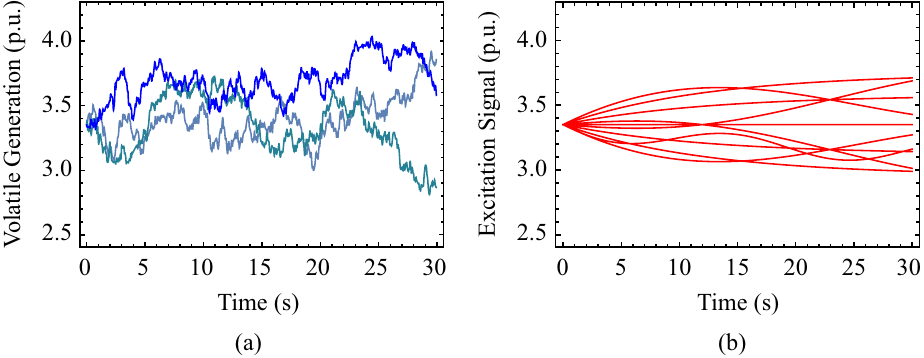}
  \caption{(a) Three paths of $P_{3}(t)$ used in the Monte Carlo simulation. (b) The excitation signals of $P_{3}(t)$ used for time-domain simulation in PSS/E with respect to the first $10$ collocation points in the proposed method, which have the greatest impact on the precision.}
  \label{fig:cps}
\end{figure}

For different orders $K$ of the KLE and error tolerances $\epsilon$ of the Smolyak algorithm, with the simulation step length $h$ set to $0.01$ s to generate excitation signals in (\ref{eq:distito}) in the MCs or in (\ref{eq:eqkl}) in the proposed method, we compare the proposed method to the MCs by finding the probability distribution and high-order moments of the RRF of the post-contingency relative rotor angle $\delta_{38\text{-}30}$ between generators 38 and 30 at $t=5$ s and the root-mean-square (RMS) value of the deviation in the system frequency $\Delta{f}$ for $t \in [20, 30]$ s. For visualization purposes, $3$ paths of the stochastic excitation $P_3(t)$ used in the MCs are shown in Fig. \ref{fig:cps}(a). The signals used for time-domain simulation in PSS/E with respect to the first $10$ collocation points in the proposed method, which have the greatest impact on the precision of the PCE, are plotted in Fig. \ref{fig:cps}(b).

Since obtaining the exact values of the expectation and moments of the RRFs is infeasible for such a complex nonlinear system, we use an MCs with a large enough sampling size $N_\mathrm{MC} = 2\times 10^6$ as a benchmark.

\begin{figure}[b]
  \centering
  \includegraphics[width=2.6in]{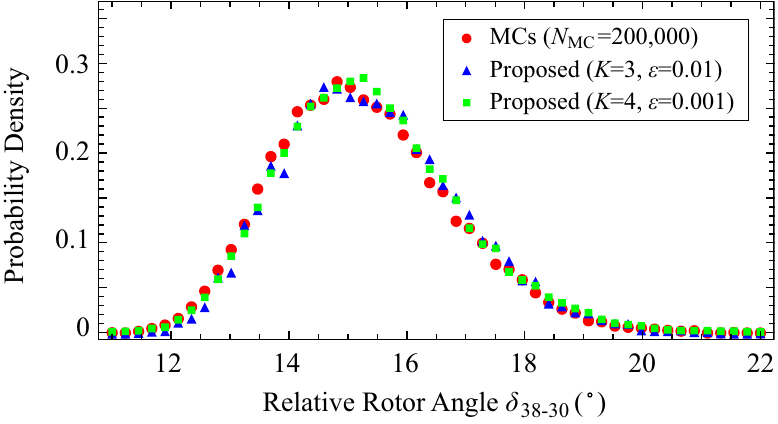}
  \caption{Probability Density of the relative rotor angles of generators $38$ and $30$ at $t=5$ s obtained by the proposed method and the MCs.}
  \label{fig:delta}\vspace{12pt}
  \centering
  \includegraphics[width=2.6in]{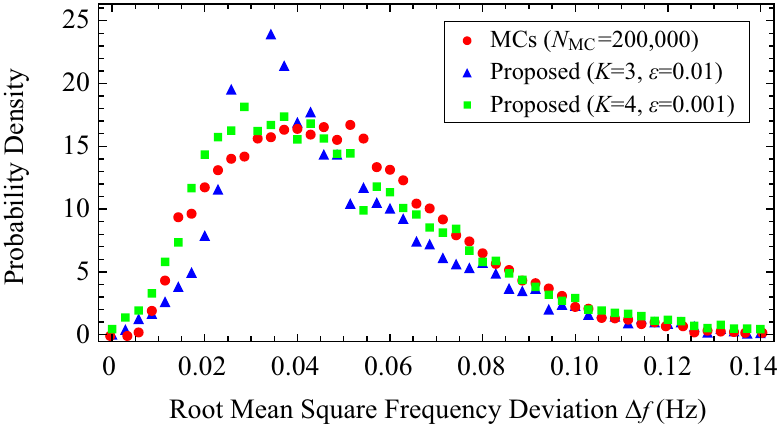}
  \caption{Probability Density of the root-mean-square deviation of the system frequency for $t\in [20,30]$ s obtained by the proposed method and the MCs.}
  \label{fig:freq}
\end{figure}

\subsubsection{Basic Results}

The probability densities of the relative rotor angle $\delta_{38\text{-}30}$ and the RMS of the frequency deviation $\Delta f$ computed using the proposed method are compared with the benchmark MCs in Figs. \ref{fig:delta} and \ref{fig:freq}, respectively. As we can see, for the relative rotor angle $\delta_{38\text{-}30}$, the proposed method produces results that are almost identical to those of the benchmark MCs with different settings of $K$ and $\epsilon$. For the RMS of the frequency deviation $\Delta f$, the result of the proposed method is not very accurate when $K = 3$ and $\epsilon=0.01$; when $K$ is increased to $4$ and $\epsilon$ is decreased to $0.001$, the obtained PDF converges to the actual one. In other words, the RRFs can be precisely extracted using the proposed method.

Furthermore, a quantitative comparison is made in Table \ref{tab:ac} by comparing the expectation and central moments of the relative angle $\delta_{38\text{-}30}$ obtained by the proposed method and the MCs with different settings. The relative error compared to the benchmark MCs is also listed. The numerical results in Table \ref{tab:ac} show that the MCs with $N_{\mathrm{MC}} = 1,000$ or $20,000$ both leads to a significant error. In contrast, the proposed method has a much better accuracy with different settings.

The computation time is also compared in Table \ref{tab:time}, where the simulation component indicates the time required for dynamic simulation in PSS/E, the method component indicates the time required for the Gaussian quadrature (\ref{eq:cm})--(\ref{eq:cmat}) and the procedure control, and the post-processing component represents the time required for the extraction of statistical information using the MCs based on the PCE.

Obviously the proposed method requires significantly less time than MCs to achieve comparable precision. This is because the key information of the uncertainty is already present in the carefully constructed excitation signals (\ref{eq:eqkl}) shown in Fig. \ref{fig:cps}(b). Therefore, a few simulations suffice to extract most of the statistical information. In fact, with $K=3$ and $\epsilon=0.01$, only $83$ simulations are needed instead of the tens of thousands required by the MCs to achieve comparable precision.

\begin{table}[tb]\scriptsize
  \renewcommand{\arraystretch}{1.25}
  \caption{Expectation, Variance, Central Moments, and Corresponding Error for the Relative Rotor Angle Obtained by the Proposed Method and the MCs In the $39$-Bus System Case} \vspace{-4.5pt}
  \label{tab:ac}
  \centering
  \begin{tabular}{p{2.4cm}p{0.75cm}p{0.75cm}p{0.785cm}p{0.785cm}p{0.85cm}}
  \hline \hline
  \hspace{24pt}\multirow{2}{*}{Method}         & \multirow{2}{*}{\hspace{6pt}Exp.}  & \multirow{2}{*}{\hspace{6pt}Var.}   & \multicolumn{3}{c}{Central Moment}        \\
  \cline{4-6}                                  &                                    &                                     & \hspace{8pt}$3$rd    & \hspace{8pt}$4$th    & \hspace{8pt}$5$th        \\
  \hline
  \hspace{0.5pt}MCs ($N_\text{MC}=2\times10^6$)
                                               & \hspace{3pt}$15.08$               & \hspace{3pt}$1.980$                & \hspace{4pt}$1.396$  & \hspace{4pt}$13.33$  & \hspace{4pt}$30.56$   \\
                                               \hline
  \multirow{2}{*}{MCs ($N_\text{MC}=20,000$)}
                                               & \hspace{3pt}$15.07$                & \hspace{3pt}$1.970$                 & \hspace{4pt}$1.573$  & \hspace{4pt}$13.96$  & \hspace{4pt}$35.89$   \\
                                               & $-0.07\%$                          & $-0.51\%$                           & $+12.7\%$            & $+5.00\%$            & $+17.4\%$             \\
                                               \hline
  \hspace{2pt}\multirow{2}{*}{MCs ($N_\text{MC}=1,000$)}
                                               & \hspace{3pt}$15.04$                & \hspace{3pt}$1.924$                 & \hspace{4pt}$1.174$  & \hspace{4pt}$11.28$  & \hspace{4pt}$19.91$   \\
                                               & $-0.28\%$                          & $-2.83\%$                           & $-15.9\%$            & $-15.4\%$            & $-34.9\%$             \\
                                               \hline
  \hspace{6pt}\multirow{2}{*}{\tabincell{c}{Proposed Method \vspace{1.5pt} \\ ($K=3$, $\epsilon=0.01$)}}
                                               & \hspace{3pt}$15.18$                & \hspace{3pt}$1.926$                 & \hspace{4pt}$1.376$  & \hspace{4pt}$12.36$  & \hspace{4pt}$26.98$   \\
                                               & $+0.64\%$                          & $-2.73\%$                           & $-1.43\%$            & $-7.27\%$            & $-11.7\%$            \\                                               \hline
  \hspace{4pt}\multirow{2}{*}{\tabincell{c}{Proposed Method \vspace{1.5pt} \\ ($K=3$, $\epsilon=0.001$)}}
                                               & \hspace{3pt}$15.17$                & \hspace{3pt}$1.982$                 & \hspace{4pt}$1.366$  & \hspace{4pt}$13.02$  & \hspace{4pt}$28.47$   \\
                                               & $+0.62\%$                          & $+0.10\%$                           & $-2.15\%$            & $-2.30\%$            & $-6.83\%$             \\
                                               \hline
  \hspace{4pt}\multirow{2}{*}{\tabincell{c}{Proposed Method \vspace{1.5pt} \\  ($K=4$, $\epsilon=0.001$)}}
                                               & \hspace{3pt}$15.17$                & \hspace{3pt}$2.045$                 & \hspace{4pt}$1.391$  & \hspace{4pt}$13.66$  & \hspace{4pt}$29.69$   \\
                                               & $+0.61\%$                          & $+3.28\%$                           & $-0.36\%$            & $+2.48\%$            & $-2.85\%$             \\

  \hline \hline
  \end{tabular}
\end{table}

\begin{table}[tb]\scriptsize
  \renewcommand{\arraystretch}{1.25}
  \caption{Computation Time of the Proposed Method and the MCs}
  \label{tab:time}
  \centering
  \begin{tabular}{cp{1cm}p{0.65cm}p{1.6cm}c}
  \hline \hline
  \multirow{2}{*}{Method}                                   & \multicolumn{3}{p{4cm}}{\hspace{28pt}Time Components (s)}       & \multirow{2}{*}{\tabincell{c}{Total \\ Time (s)}}       \\
  \cline{2-4}                                               & Simulation               & Method        & Post-Processing    &                            \\
  \hline
  MCs ($N_\text{MC}=20,000$) \vspace{1.5pt}                 &\hspace{2pt}$2207.97$     &\hspace{7pt}$-$           &\hspace{15pt}$0.55$        & $2208.52$                  \\
  MCs ($N_\text{MC}=1,000$) \vspace{2pt}                    &\hspace{3.5pt}$115.58$    &\hspace{7pt}$-$           &\hspace{15pt}$0.15$        & $115.73$   \\
  \tabincell{c}{Proposed Method \vspace{1pt}\\ ($K=3$, $\epsilon=0.01$) } \vspace{3pt}
                                                            &\hspace{8pt}$9.97$        &\hspace{3pt}$0.05$        &\hspace{15pt}$1.23$        & $11.25$ \\
  \tabincell{c}{Proposed Method  \vspace{1pt}\\ ($K=3$, $\epsilon=0.001$)} \vspace{3pt}
                                                            &\hspace{6pt}$37.85$       &\hspace{3pt}$0.17$        &\hspace{15pt}$3.03$        & $41.05$ \\
  \tabincell{c}{Proposed Method  \vspace{1pt}\\ ($K=4$, $\epsilon=0.001$)} \vspace{0.3pt}
                                                            &\hspace{6pt}$55.87$       &\hspace{3pt}$0.26$        &\hspace{15pt}$4.88$        & $61.01$ \\
  \hline \hline
  \end{tabular}
\end{table}

\subsubsection{Comprehensive Comparison}

\begin{figure}[t]
  \centering
  \includegraphics[width=2.74in]{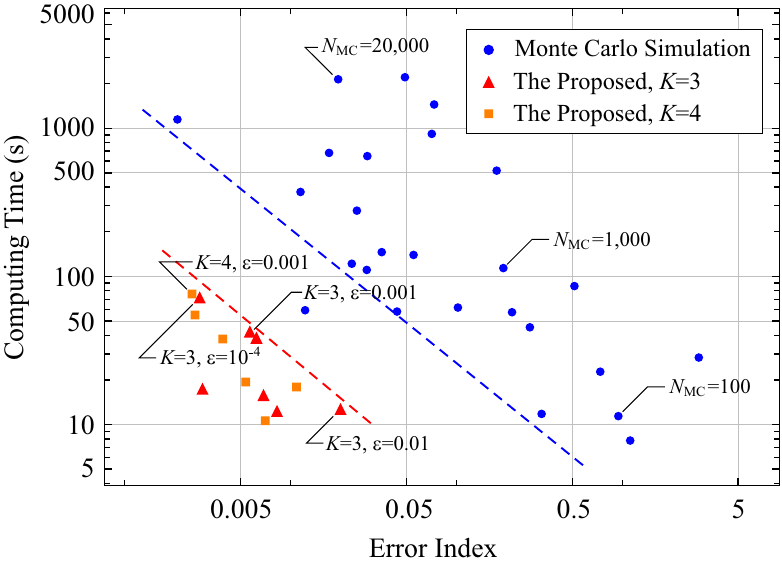}
  \caption{Computation time versus error index for the proposed method and the Monte Carlo simulation with different settings for the $39$-bus system.}
  \label{fig:terr}
\end{figure}

Since both the proposed method and the MCs have adjustable precision settings, i.e., $K$ and $\epsilon$ in the proposed method and $N_{\mathrm{MC}}$ in the MCs, to compromise between precision and the computational burden, we repeatedly apply them with different settings to make a comprehensive comparison. To quantify the precision, we define an error index as the sum of the squares of the relative errors of the expectation and the first $5$ central moments, as:
\begin{align}
    err = & \left[({E - \hat{E}})/{E}\right]^2  + \sum_{j=2}^5 \left[{ (M^{(j)} - \hat{{M}}^{(j)})}/{M^{(j)}} \right]^2 \label{eq:err}
\end{align}
\noindent
where $E$ and $\hat{E}$ represent the benchmark and obtained expectation of the RRF; $M^{(j)}$ and $\hat{M}^{(j)}$ represent those of the $j$th central moment.

Fig. \ref{fig:terr} shows the computation time versus the precision of the proposed method and the MCs in logarithmic coordinates. Although both methods exhibit negative correlations between the error index and the computation time, the results of the proposed method are more than one to two orders of magnitude better than those of the MCs. Moreover, the results of the MCs vary over a very wide range for a given setting, which is inherent in the nature of random sampling. In contrast, the proposed method performs much more consistently, as only deterministic simulations are required to find the PCE.

\begin{figure}[t]
  \centering
  \includegraphics[width=3.48in]{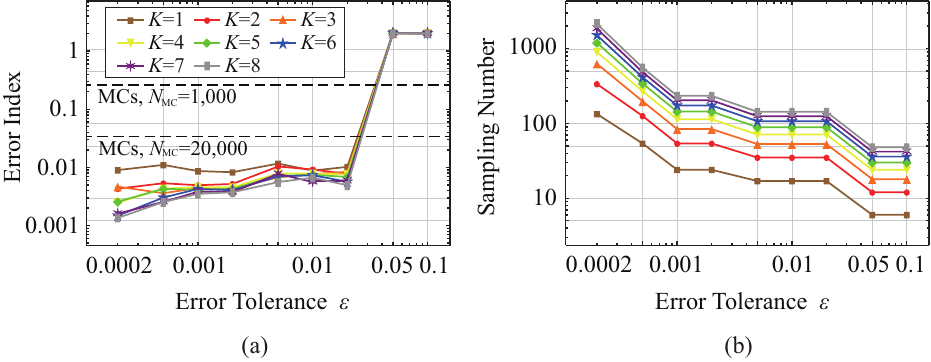}
  \caption{Error index and sampling number of the proposed method with different $\epsilon$ and $K$ in evaluating the relative rotor angle $\delta_{38\text{-}30}$ in the IEEE 39-bus system. (a) Error index. (b) Sampling number. Error indices of the MCs are also plotted for comparison.}
  \label{fig:order39}
\end{figure}

\begin{figure}[b]
  \centering
  \includegraphics[width=2.6in]{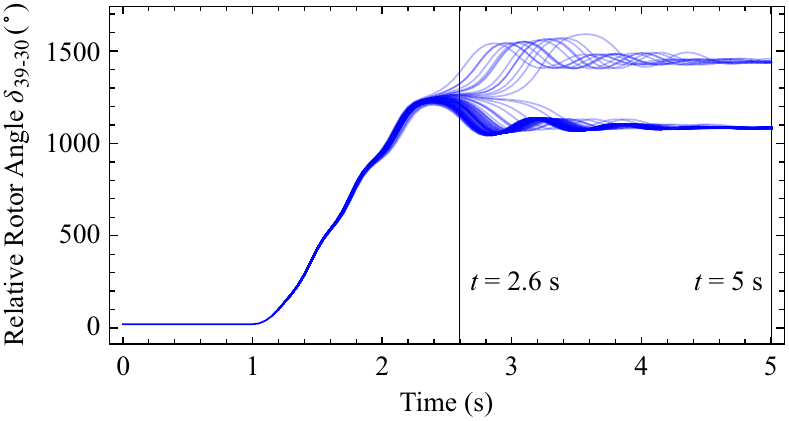}
  \caption{$100$ random paths of the relative rotor angle between generators $39$ and $30$ for $t\in[0,5]$ s when instability occurs.}
  \label{fig:traj}\vspace{12pt}
  \includegraphics[width=3.47in]{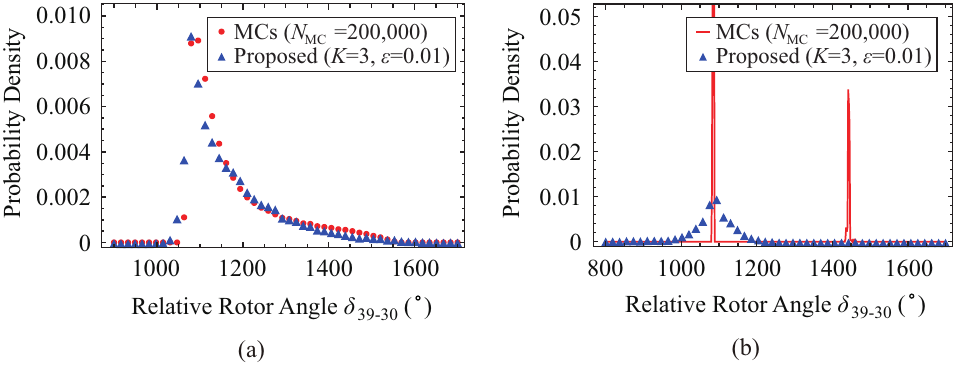}
  \caption{Probability Density of the relative rotor angles of generators 39 and 30 obtained by the proposed method and MCs. (a) $t = 2.6$ s. (b) $t = 5$ s.}
  \label{fig:instability}
\end{figure}

\subsubsection{Numerical Convergence Analysis}
\label{sec:nconv}

The impact of $K$ and $\epsilon$ on the accuracy of the proposed method is investigated. Figs. \ref{fig:order39}(a) and \ref{fig:order39}(b) respectively present the error index defined by (\ref{eq:err}) and sampling number of the proposed method with different $\epsilon$ and $K$ in evaluating the relative rotor angle $\delta_{38\text{-}30}$. We can observe that as $\epsilon$ is decreased below $0.05$, the error index rapidly drop below $0.01$. However, further decreasing $\epsilon$ only leads to a slower decrease of the error index. This phenomenon can be attributed to the error introduced by truncating the KL expansion (\ref{eq:kl}), the random error of using MCs to extract statistical information from the PCE, and possibly numerical error of the simulation software. In contrast, Fig. \ref{fig:order39}(b) shows that as $\epsilon$ is further decreased, the computation burden represented by the sampling number increases significantly. From these figures, we can see that as discussed in Sections \ref{sec:kl} and \ref{sec:pcm}, choosing $\epsilon$ from $0.01$ to $0.001$ and $K$ from $3$ to $6$ achieves a satisfactory balance between the precision and computational burden.

Note that higher $K$ represents higher frequencies of the alternating components in the disturbance signal. Hence, the convergence of the error index in terms of $K$ also reflects the fact that treating the random disturbances as stochastic processes instead of random parameters is necessary.

\subsubsection{Case of Instability} \label{sec:instability}
The proposed method is also tested when instability occurs. To induce instability, a three-phase grounding fault is imposed on bus 39 at t = 1 s and lasts for 0.2 s until it is cleared by tripping line 39-9. Other settings are the same as the previous case. The post-fault relative rotor angle between generators $39$ and $30$, i.e., $\delta_{39\text{-}30}$, is investigated.

For easy understanding, we plot $100$ random trajectories of $\delta_{39\text{-}30}(t)$ for $t\in[0,5]$s in Fig. \ref{fig:traj}. The proposed method with $K=3$ and $\epsilon=0.01$ is used to find the probability density of $\delta_{39\text{-}30}$ at $t=2.6$ s and $5$ s, respectively. The results compared to the MCs are plotted in Figs. \ref{fig:instability}.

For $t=2.6$ s, albeit the nonlinearity of an instable power system, the proposed method precisely captures the probability density of the divergent rotor angle, where only $153$ simulations are performed. The first $3$ orders of moments obtained are $1172.7$, $1.389\times10^6$, and $1.660\times10^9$, which are almost identical to the results of $1171.8$, $1.383\times10^6$, and $1.646\times10^9$ obtained by MCs with $N_{\mathrm{MC}}=2\times10^6$.

As of $t=5$ s, from Fig. \ref{fig:instability}(b), we can see the probability density becomes two spikes. The height of the left spike is $0.298$ and the right one is $0.032$. Unfortunately, the proposed method only captures left one. This is because the RRF is almost no more a continuous function of $\bm{\zeta}$. From (\ref{eq:appconv}), we can see that the convergence rate of PCE deteriorates to the worst case. Hence, the proposed method only captures the dominant spike with a limited approximation order. If one wants to capture both two spikes, a very high approximation order is needed, which is impractical in application. This is currently a limitation of the proposed method, as noted in Section \ref{sec:conv}.

\subsection{IEEE 118-Bus System Case Study}
\label{sec:118}

\begin{figure}[tb]
  \centering
  \includegraphics[width=3.5in]{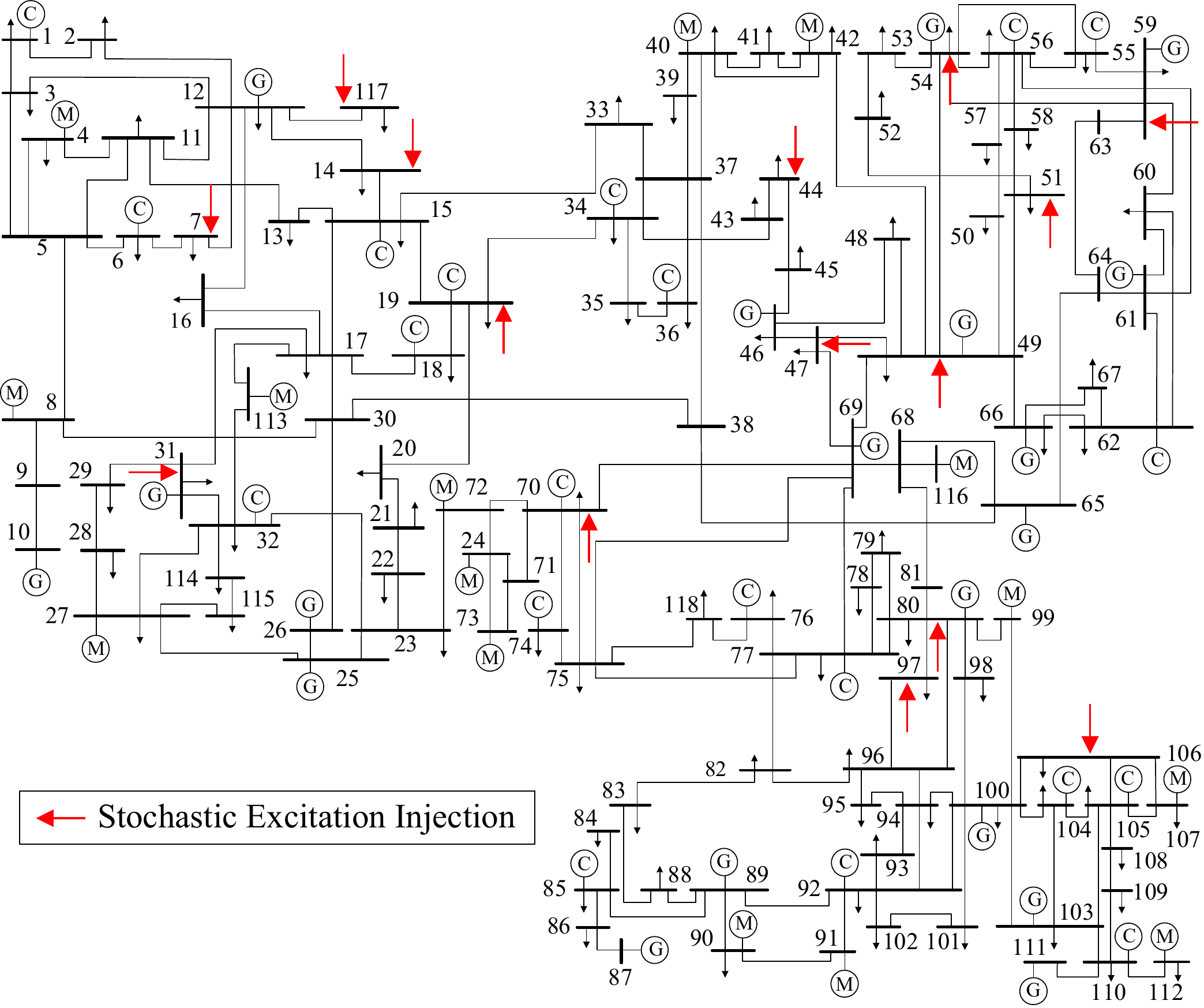}
  \caption{Diagram of the 118-bus system. The injection points of the $15$ stochastic excitations are labeled.}
  \label{fig:118sys}
\end{figure}

\subsubsection{Setting}

To further verify the proposed method in larger systems with high-dimensional stochastic excitation inputs, we test it on the IEEE 118-bus system \cite{118sys}. Fifteen stochastic excitations defined by $5$ independents sets of (\ref{eq:case}) are connected as power injections to buses $54$, $59$, $80$, $7$, $14$, $117$, $19$, $31$, $46$, $47$, $49$, $70$, $44$, $51$, and $97$, as shown in Fig. \ref{fig:118sys}.
A three-phase grounding is applied to bus $49$ at $t=1$ s and lasts for $0.2$ s until it is cleared by tripping line $49\text{-}69$. The uncertainty in the relative rotor angle $\delta_{69\text{-}49}$ between generators $69$ and $49$ at $t=5$ s is quantified.

\subsubsection{Basic Results}

Let $K=3$ and $\epsilon=0.01$, the proposed method takes $293.8$ s to find the expectation and moments of $\delta_{69\text{-}49}$. In the proposed method, $941$ simulations are performed on PSS/E to evaluate the collocation points. Meanwhile, the MCs with $N_\mathrm{MC} = 2,000$ takes $597.1$ s, and the MCs with $N_\mathrm{MC} = 20,000$ takes $6113.6$ s, as listed in Table \ref{tab:time118}. Moreover, using result of the MCs with $N_\mathrm{MC} = 2\times10^6$ as a benchmark, the error index defined by (\ref{eq:err}) for the proposed method is $3.26\times 10^{-3}$, which is substantially better than the value of $0.227$ for the MCs with $N_\mathrm{MC} = 2,000$, and is also much better than the value of $2.66\times 10^{-2}$ for the the MCs with $N_\mathrm{MC} = 20,000$. Detailed comparison of the calculated statistical information is given in Table \ref{tab:118ac}. Obviously, the proposed method significantly outperforms the MCs in this larger power system with higher-dimensional stochastic inputs.

\begin{table}[tb]\scriptsize
  \renewcommand{\arraystretch}{1.35}
  \caption{Computation Time of the Proposed Method and the MCs in the IEEE 118-Bus System Case}
  \label{tab:time118}
  \centering
  \begin{tabular}{cp{1cm}p{0.65cm}p{1.6cm}c}
  \hline \hline
  \multirow{2}{*}{Method}                                   & \multicolumn{3}{p{4cm}}{\hspace{28pt}Time Components (s)}       & \multirow{2}{*}{\tabincell{c}{Total \\ Time (s)}}       \\
  \cline{2-4}                                               & Simulation               & Method        & Post-Processing    &                            \\
  \hline
  MCs ($N_\text{MC}=20,000$) \vspace{1.6pt}                 &\hspace{4.5pt}$6113.0$     &\hspace{7pt}$-$           &\hspace{18pt}$0.5$        & $6113.6$                  \\
  MCs ($N_\text{MC}=2,000$) \vspace{2pt}                    &\hspace{6.5pt}$596.9$    &\hspace{7pt}$-$           &\hspace{18pt}$0.2$        & $597.1$   \\
  \tabincell{c}{Proposed Method \vspace{1pt}\\ ($K=3$, $\epsilon=0.01$) } \vspace{0.3pt}
                                                            &\hspace{6.5pt}$285.2$     &\hspace{5.5pt}$1.1$        &\hspace{18pt}$7.5$        & $293.8$ \\  \hline \hline
  \end{tabular}
\end{table}

\begin{table}[tb]\scriptsize
  \renewcommand{\arraystretch}{1.35}
  \caption{Expectation, Variance, Central Moments, and Corresponding Error for the Relative Rotor Angle Obtained by the Proposed Method and the MCs In the $118$-Bus System Case}
  \label{tab:118ac}
  \centering
  \begin{tabular}{p{2.4cm}p{0.75cm}p{0.75cm}p{0.785cm}p{0.785cm}p{0.85cm}}
  \hline \hline
  \hspace{24pt}\multirow{2}{*}{Method}         & \multirow{2}{*}{\hspace{6pt}Exp.}  & \multirow{2}{*}{\hspace{6pt}Var.}   & \multicolumn{3}{c}{Central Moment}        \\
  \cline{4-6}                                  &                                    &                                     & \hspace{8pt}$3$rd    & \hspace{8pt}$4$th    & \hspace{8pt}$5$th        \\
  \hline
  \hspace{0.5pt}MCs ($N_\text{MC}=2\times10^6$)
                                               & \hspace{3pt}$17.24$                & \hspace{3pt}$0.749$                 &  $-0.137$            & \hspace{4pt}$1.700$  & $-1.012$   \\
                                               \hline
  \multirow{2}{*}{MCs ($N_\text{MC}=20,000$)}
                                               & \hspace{3pt}$17.22$                & \hspace{3pt}$0.740$                 & $-0.176$             & \hspace{4pt}$1.727$  & $-1.157$   \\
                                               & $-0.10\%$                          & $-1.12\%$                           & $+7.50\%$            & $+1.60\%$            & $+14.35\%$             \\
                                               \hline
  \hspace{2pt}\multirow{2}{*}{MCs ($N_\text{MC}=2,000$)}
                                               & \hspace{3pt}$17.23$                & \hspace{3pt}$0.775$                 & $-0.180$             & \hspace{4pt}$1.950$  & $-1.341$   \\
                                               & $+0.04\%$                          & $+3.46\%$                           & $+31.42\%$           & $+14.72\%$            & $+32.54\%$             \\
                                               \hline
  \hspace{6pt}\multirow{2}{*}{\tabincell{c}{Proposed Method \vspace{1.5pt} \\ ($K=3$, $\epsilon=0.01$)}}
                                               & \hspace{3pt}$17.20$                & \hspace{3pt}$0.726$                 & $-0.134$             & \hspace{4pt}$1.631$  &$-1.033$   \\
                                               & $-0.22\%$                          & $-2.96\%$                           & $-1.68\%$            & $-4.06\%$            & $+2.10\%$            \\                                               \hline\hline
  \end{tabular}
\end{table}

\begin{table}[tb]\scriptsize
  \renewcommand{\arraystretch}{1.35}
  \caption{Computation Time of the Proposed Method with Different Dimension of Stochastic Excitations in the IEEE 118-Bus System Case}
  \label{tab:times118}
  \centering
  \begin{tabular}{ccccc}
  \hline \hline
  \multirow{2}{*}{\tabincell{c}{Random \\ Inputs}}    & \multicolumn{3}{c}{Time Components (s)}       & \multirow{2}{*}{\tabincell{c}{Total \\ Time (s)}}       \\
  \cline{2-4}                       & Simulation               & Method       & Post-Processing    &                            \\
  \hline
  $3$                               &$22.9$        &$0.1$     &$1.3$        & $24.3$ \\
  $6$                               &$71.2$        &$0.2$     &$2.8$        & $74.1$ \\
  $9$                               &$138.1$       &$0.6$     &$4.3$        & $143.0$ \\
  $12$                              &$199.7$       &$0.8$     &$4.8$        & $205.4$ \\
  $15$                              &$285.2$       &$1.1$     &$7.5$        & $293.8$ \\
  \hline \hline
  \end{tabular}
\end{table}

In addition, we perform a numerical convergence analysis similar to Section \ref{sec:nconv}. The result is substantially the same, showing choosing $\epsilon$ from $0.01$ to $0.001$ and $K$ from $3$ to $6$ achieves a satisfactory balance between precision and computation burden. In order not to make this paper too lengthy, detailed result will not be given here.

\subsubsection{Discussions on Computation Cost versus the Dimension of Stochastic Excitation Inputs} Finally, we investigate the relation between the computation time consumed by be proposed method and the dimension of the stochastic excitations. With the dimension ranging between $3$ and $15$, the computation time is shown in Table \ref{tab:times118}. We can see that the computation time of the proposed method is almost linear with respect to the dimension of random inputs, which means the proposed method is able to handle tens of random inputs. However, as the inherent limitation of a PCE-based method  \cite{sun2018probabilistic,xu2019probabilistic}, the proposed method cannot accommodate very large number of random inputs. To deal with this situation, it could be possible to use dimension reduction techniques to reduce the number of independent random driving forces of the excitations (the standard Wiener processes $\bm{W}_t$ in this case), which may further improve the efficiency of the proposed method.

\section{Conclusions}

A nonintrusive method for quantifying the uncertainty in dynamic power systems based on SDE and PCE is proposed in this paper. The proposed method exhibits high precision and efficiency compared to the Monte Carlo simulation. The method runs on commercial simulation software such as PSS/E, which ensures ease of use for power utilities.

Currently, discrete switching events or extreme nonlinearity may negatively impact the convergence rate of the proposed method, which may lead to inaccurate results in certain situations. To alleviate such limitation, leveraging multi-element PC technique in the framework of the proposed method is worth investigating in the future study.

Extending the proposed method from the It\^{o} process driven by Wiener processes relevant to Hermite PC to other random processes related to the generalized polynomial chaos (gPC) to handle more complicated stochastic phenomena in power systems is also one of the promising directions for future research.

\appendices

\section{The Fokker-Planck Equation}
\label{sec:fp}
\setcounter{equation}{0}
\renewcommand{\theequation}{\ref{sec:fp}\arabic{equation}}

The probability density of the It\^{o} process evolving with time admits the following partial differential equation, known as the \emph{Fokker-Planck equation}, formulated as
\begin{align}
    \frac{\partial p\left( \bm{\xi}_t, t \right)}{\partial t}  = & - \sum_{j=1}^{m} \frac{\partial}{\partial \xi_{j}} \left[ \mu_j \left( \bm{\xi}_t, t \right) p\left( \bm{\xi}_t, t \right)\right] \nonumber \\
    & + \sum_{j=1}^{m} \sum_{k=1}^{m} \frac{\partial^2}{\partial \xi_{j} \partial \xi_{k}} \left[ D_{jk}\left( \bm{\xi}_t, t \right) p\left( \bm{\xi}_t, t \right)\right], \label{eq:fp}
\end{align}

\noindent
where $p(\cdot)$ is the probability density function (PDF) of $\bm{\xi}_t$; $\xi_j$ is the $j$th entry of $\bm{\xi}_t$; $\bm{D}(\cdot) = \bm{\sigma}(\cdot) \bm{\sigma}^{\mathrm{T}}(\cdot) / 2$; and $D_{jk}(\cdot)$ is the $k$th entry in the $j$th row of $\bm{D}(\cdot)$.
By setting the left-hand side to zero, the stationary state is represented.

\section{Method for Identifying the It\^{o} Process Model From Real Data}
\label{sec:data}
\setcounter{equation}{0}
\renewcommand{\theequation}{\ref{sec:data}\arabic{equation}}

Suppose a set of recorded data of renewable generations with sampling interval $h$, denoted by $\{ \tilde{\bm{\xi}}_0, \tilde{\bm{\xi}}_h,\tilde{\bm{\xi}}_{2h},\ldots,\tilde{\bm{\xi}}_T \}$. Construct the drift and diffusion terms $\bm{\mu} \left(\bm{\xi}_t,t ; \bm{q} \right)$ and $\bm{\sigma} \left(\bm{\xi}_t, t; \bm{q} \right)$ in (\ref{eq:ito}) as simple functions of $\bm{\xi}_t$, such as polynomials with parameters $\bm{q}$ to be identified, such that the likelihood of the following logarithmic conditional probability is maximized:
\begin{align}
    \max_{\bm{q}} L = \log \Pr \left[ \tilde{\bm{\xi}}_h,\tilde{\bm{\xi}}_{2h},\ldots,\tilde{\bm{\xi}}_T | \tilde{\bm{\xi}}_{0} \right].  \label{eq:likelihood}
\end{align}

By the independent incremental property of the It\^{o} process \cite{Pardoux2014Stochastic}, the conditional probability in (\ref{eq:likelihood}) can be rewritten as
\begin{align}
    L = - \sum\nolimits_{j=1}^{T/h} \log \Pr \left[ \tilde{\bm{\xi}}_{jh} | \tilde{\bm{\xi}}_{(j-1)h} \right]. \label{eq:like}
\end{align}

Considering the discrete form (\ref{eq:distito}) of the It\^{o} process (\ref{eq:ito}), and considering that the sampling interval $h$ is short, we obtain
\begin{align}
   \bm{\xi}_{t+h} \sim \mathcal{N} \left( \bm{\xi}_t + h  \bm{\mu}_t, h \bm{\sigma}^\mathrm{T}_t  \bm{\sigma}_t \right).
\end{align}
\noindent
where ${\bm{\mu}}_t$ and ${\bm{\sigma}}_t$ represent $\bm{\mu}( {\bm{\xi}}_{t},t; \bm{q})$ and $\bm{\sigma}( {\bm{\xi}}_{t},t; \bm{q})$.

Therefore, the conditional probability in (\ref{eq:like}) is
\begin{align}
   \Pr & \Big[  \tilde{ \bm{\xi}}_{t+h}  | \tilde{\bm{\xi}}_{t}  \Big]   = \dfrac{1}{\sqrt{\left( 2 \pi h \right)^m \det\left( \tilde{\bm{\sigma}}_t\tilde{\bm{\sigma}}_t^\mathrm{T}\right) }}  \label{eq:reprob} \\
   \times   & \exp   \left\{ -\frac{1}{2}\left[  \Delta\tilde{\bm{\xi}}_{t} - h\tilde{\bm{\mu}}_t \right]^\mathrm{T} \left( \tilde{\bm{\sigma}}_t\tilde{\bm{\sigma}}_t^\mathrm{T}\right)^{-1} \left[ \Delta\tilde{\bm{\xi}}_{t} - h\tilde{\bm{\mu}}_t \right] \right\}, \nonumber
\end{align}
\noindent
where $\Delta \tilde{\bm{\xi}}_{t} = \tilde{\bm{\xi}}_{t+h} - \tilde{\bm{\xi}}_{t}$ represents the change in the recorded data over a sampling interval; $\tilde{\bm{\mu}}_t$ and $\tilde{\bm{\sigma}}_t$ represent $\bm{\mu}( \tilde{\bm{\xi}}_{t},t; \bm{q})$ and $\bm{\sigma}( \tilde{\bm{\xi}}_{t},t; \bm{q})$, respectively.

Substituting (\ref{eq:reprob}) into (\ref{eq:like}), letting $\tilde{\bm{D}}_{jh} =h \tilde{\bm{\sigma}}_{jh} \tilde{\bm{\sigma}}_{jh}^{\mathrm{T}}/2$, and neglecting the constant terms in the logarithmic function yields
\begin{align}
  \min_{\bm{q}} L'  = & \frac{1}{4} \sum_{j=0}^{\frac{T}{h}-1} \left[  \Delta\tilde{\bm{\xi}}_{jh} - h\tilde{\bm{\mu}}_{jh} \right]^\mathrm{T} \tilde{\bm{D}}_{jh}^{-1} \left[  \Delta\tilde{\bm{\xi}}_{jh} - h\tilde{\bm{\mu}}_{jh} \right] \nonumber \\
  + & \frac{1}{2}  \sum_{j=0}^{\frac{T}{h}-1} \log \det \left(\tilde{\bm{D}}_{jh}\right). \label{eq:ident}
\end{align}

Model (\ref{eq:ident}) is an unconstrained programming. Although not necessarily convex, generally it can be solved using common methods such as the gradient descent. In case of non-convexity, heuristic global optimization methods such as the particle swarm optimization (PSO) can also be applied. As (\ref{eq:ident}) is solved, the optimal parameters of the It\^{o} process model are found.

\section{Proof of Proposition \ref{prop:2}}
\label{sec:proof}
\setcounter{equation}{0}
\renewcommand{\theequation}{\ref{sec:proof}\arabic{equation}}

Based on the triangular inequality, we have
\begin{align}
 & \| \omega \left( \{ \bm{\xi}_\tau \}_{\tau\in[0,t]}\right) - \hat{\omega}\left( \bm{\zeta} \right) \|    \label{eq:tri} \\
 & \hspace{10pt} \le    \| \omega\left( \{ \bm{\xi}_\tau \}_{\tau\in[0,t]}\right)  - \omega^*\left( \bm{\zeta} \right)  \| + \|  \omega^*\left( \bm{\zeta} \right) -   \hat{\omega}\left( \bm{\zeta} \right)  \|. \nonumber
\end{align}

As of the first term on the right-hand side of (\ref{eq:tri}), jointly using Assumption \ref{assu:1} and Proposition \ref{prop:1}, we have
\begin{align}
  & \| \omega\left( \{ \bm{\xi}_\tau \}_{\tau\in[0,t]}\right)  - \omega^*\left( \bm{\zeta} \right)  \|  \label{eq:bound} \\
  & \hspace{22pt} \le  L \| \{\bm{\xi}_\tau  -  \bm{\xi}^*_\tau\left( \bm{\zeta} \right) \}_{\tau\in[0,t]} \| \rightarrow 0, \text{as } K \rightarrow \infty. \nonumber
\end{align}

As of the second term on the right-hand side of (\ref{eq:tri}), polynomial approximation theory \cite{Fletcher1984,Xiu2010} indicates that there exists a constant $C$ independent of $N$ and $p \in \mathbb{N}^{+}$ such that
\begin{align}
  & \| \omega^*\left( \bm{\zeta} \right) - \hat{\omega}\left( \bm{\zeta} \right)  \|  \label{eq:appconv}\\
  & \le \Bigg\{ \begin{array}{ll}
  C e^{-\alpha N} \| \omega^*\left( \bm{\zeta} \right) \|,  &\hspace{-6pt}\omega^* \text{ is analytical}, \\
  C N^{-p} \| \omega^*\left( \bm{\zeta} \right) \|_{H^p_w}, &\hspace{-6pt}\omega^* \in H^p_w, \\
  C N^{\frac{1}{2}} \| \omega^*\left( \bm{\zeta} \right) \|, &\hspace{-6pt}\omega^*\text{ has finite discontinuity points},
  \end{array}  \nonumber
\end{align}
where $N$ is the order of PCE $\hat{\omega}\left( \bm{\zeta} \right)$ with a complete polynomial basis; ${H^p_w}$ is a Sobolev space of order $p$, as
\begin{align}
    H^k_w \triangleq \Big\{ v(\bm{\zeta}): \frac{d^{|\bm{j}|}v}{d\zeta_1^{j_1}d\zeta_2^{j_2}\cdots d\zeta_M^{j_M}} \in L^2_w, |\bm{j}| \le k \Big\}, \label{eq:sobo}
\end{align}
\noindent
and $\| \cdot \|_{H^k_w}$ is the norm for $v \in {H^k_w}$, defined as
\begin{align}
   \Vert v \Vert_{H^k_w} \triangleq \Bigg[ \sum_{|\bm{j}| \le k} \bigg\Vert \frac{d^{|\bm{j}|}v}{d\zeta_1^{j_1}d\zeta_2^{j_2}\cdots d\zeta_M^{j_M}} \bigg\Vert \Bigg]^{1/2}.
\end{align}

From (\ref{eq:appconv}) we can see that the convergence rate of the PCE depends on the smoothness of $\omega^*\left( \cdot \right)$. When $\omega^*\left( \cdot \right)$ is analytical (infinitely differentiable), the approximation has and exponential convergence rate. Moreover, providing that Assumptions \ref{assu:2} and $\ref{assu:3}$ hold, (\ref{eq:appconv}) converges to $0$ as $N\rightarrow \infty$.

The Smolyak adaptive sparse collocation algorithm ensures the equivalence between $\epsilon\rightarrow 0$ and  $N\rightarrow \infty$ \cite{CONRAD2013ADAPTIVE}. Due to the proof is rather tedious, we will not present it here, as it can be found in related literature such as \cite{CONRAD2013ADAPTIVE}. Thus, we obtain
\begin{align}
  \| \omega^*\left( \bm{\zeta} \right) - \hat{\omega}\left( \bm{\zeta} \right)  \| \rightarrow 0,\  \text{as } \epsilon \rightarrow 0,  \label{eq:appconv1}
\end{align}

Combining (\ref{eq:tri}), (\ref{eq:bound}), and (\ref{eq:appconv1}) concludes this proof.

\bibliographystyle{IEEEtran}
\bibliography{IEEEabrv,SAPC}


\begin{IEEEbiography}[{\includegraphics[width=1in,height=1.25in,clip,keepaspectratio]{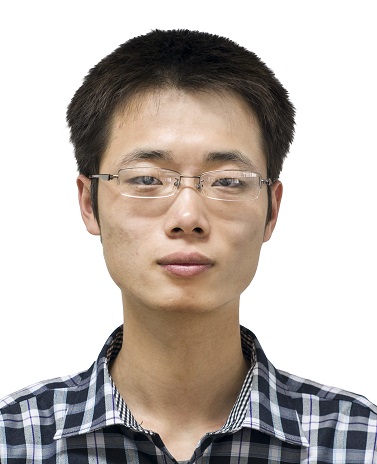}}]
{Yiwei Qiu}(S'17-M'18)
was born in Taizhou, Zhejiang, China in 1991. He received the B.S. and Ph.D. degrees from Zhejiang University, Hangzhou, Zhejiang, China, respectively in 2013 and 2018, both in Engineering. \\
\indent He is currently a Post-Doctoral Research Fellow at Department of Electrical Engineering, Tsinghua University, Beijing, China. He was a Visiting Research Associate at Department of Electrical Engineering, Southern Methodist University, TX, USA, from 2017 to 2018. His research interests include power system dynamics and control, and integration of renewable energy.
\end{IEEEbiography}

\begin{IEEEbiography}  [{\includegraphics[width=1in,height=1.25in,clip,keepaspectratio]{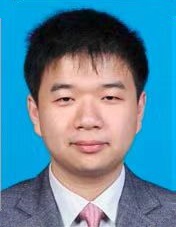}}]
{Jin Lin}(S'11-M'12)
was born in 1985. He received the B.S. and Ph.D. degrees in electrical engineering from Tsinghua University, Beijing, China, in 2007 and 2012, respectively. \\
\indent He is currently an Associate Professor with the Department of Electrical Engineering, Tsinghua University, Beijing, China. His research interests are in grid integration technology of renewable energy and dynamic power systems.
\end{IEEEbiography}

\begin{IEEEbiography} [{\includegraphics[width=1in,height=1.25in,clip,keepaspectratio]{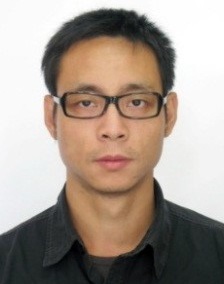}}]
{Feng Liu}(M'10-SM'18)
received the B.Sc. and Ph.D. degrees in electrical engineering from Tsinghua University, Beijing, China, in 1999 and 2004, respectively. \\
\indent Dr. Liu is currently an Associate Professor of Tsinghua University. From 2015 to 2016, he was a visiting associate at California Institute of Technology, CA, USA. His research interests include
stability analysis, optimal control, robust dispatch and game theory based decision making in energy and power systems. He is the author/coauthor of more than 200 peer-reviewed technical papers and two books, and holds more than 20 issued/pending patents. He was a guest editor of IEEE Transactions on Energy Conversion.
\end{IEEEbiography}

\begin{IEEEbiography} [{\includegraphics[width=1in,height=1.25in,clip,keepaspectratio]{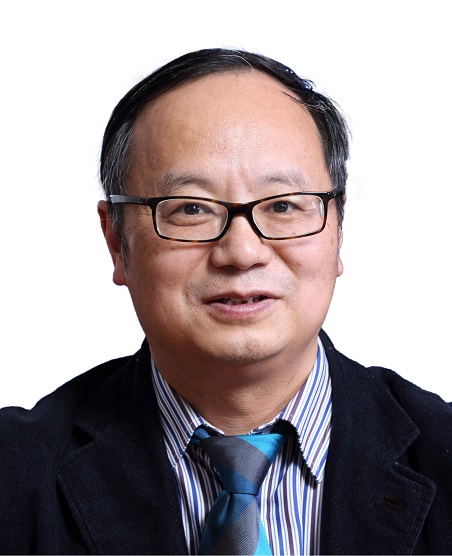}}]
{Yonghua Song} (M'90$-$SM'94$-$F'08) received the B.Eng. and Ph.D. degrees from Chengdu University of Science and Technology, Chengdu, China, and China Electric Power Research Institute, Beijing, China, in 1984 and 1989, respectively.
He is a Professor with Department of Electrical and Computer Engineering, The University of Macau, and also an Adjunct Professor at Department of Electrical Engineering, Tsinghua University.
His current research interests include smart grid, electricity economics, and operation and control of power systems.\\
\indent Prof. Song was a recipient of the D.Sc. Award from Brunel University in 2002 for his original achievements in power systems research. He was elected Vice-President of the Chinese Society for Electrical Engineering (CSEE) and appointed Chairman of the International Affairs Committee of the CSEE in 2009. He was elected as a Fellow of the Royal Academy of Engineering, U.K. in 2004.
\end{IEEEbiography}

\vfill

\end{document}